\documentclass[a4paper]{amsart}

\usepackage{amsmath,amsfonts,amssymb,amsthm} 
\usepackage{mathrsfs} 
\usepackage[shortlabels]{enumitem} 
\usepackage{caption} 
\usepackage[T1]{fontenc} 

\usepackage{mathtools} 
\mathtoolsset{centercolon} 

\usepackage{tikz} 
\usetikzlibrary{cd} 

\usepackage{pdflscape,multirow}

\usepackage{hyperref} 

\usepackage{todonotes}

\theoremstyle{definition}
\newtheorem{definition}{Definition}[section]
\newtheorem{assumption}[definition]{Assumption}

\theoremstyle{plain}
\newtheorem{lemma}[definition]{Lemma}
\newtheorem{theorem}[definition]{Theorem}
\newtheorem{proposition}[definition]{Proposition}
\newtheorem{corollary}[definition]{Corollary}
\newtheorem{conjecture}[definition]{Conjecture}

\theoremstyle{remark}
\newtheorem{remark}[definition]{Remark}
\newtheorem{example}[definition]{Example}

\newtheorem{method}[definition]{Method}

\newcommand{\bA}{\mathbb{A}}

\newcommand{\bC}{\mathbb{C}}

\newcommand{\bF}{\mathbb{F}}

\newcommand{\bP}{\mathbb{P}}
\newcommand{\bQ}{\mathbb{Q}}
\newcommand{\bR}{\mathbb{R}}

\newcommand{\bZ}{\mathbb{Z}}
\newcommand{\WP}{\mathbb{W}\mathbb{P}}

\newcommand{\calL}{\mathcal{L}}

\newcommand{\calX}{\mathcal{X}}

\newcommand{\rmI}{\mathrm{I}}

\newcommand{\rmK}{\mathrm{K}}

\newcommand{\rmV}{\mathrm{V}}

\DeclareMathOperator{\MW}{MW}

\DeclareMathOperator{\NS}{NS}

\DeclareMathOperator{\rank}{rank}

\DeclareMathOperator{\SN}{SN}

\DeclareMathOperator{\Trv}{Trv}

\numberwithin{figure}{section}
\numberwithin{table}{section}

\begin{document}

\title[Computational Tools for Families of Elliptic Surfaces]{Computational Tools for Families of Elliptic Surfaces: N\'{e}ron-Severi and Mirror Symmetry}

\author[C. F. Doran]{Charles F. Doran}
\address{Department of Mathematical and Statistical Sciences, University of Alberta, University Commons 5-140, Edmonton, Alberta, T6G 2N8, Canada.}
\email{charles.doran@ualberta.ca}
\thanks{The research of CFD was supported by the Natural Sciences and Engineering Research Council of Canada through the NSERC Discovery Grant ``Fibration and Degeneration in Calabi-Yau Geometry'' (RGPIN-2024-05855), and by the NSERC Alliance-Alberta Innovates Advance Program (252608224). CFD also acknowledges support from the Center of Mathematical Sciences and Applications at Harvard University and the Distinguished Visiting Professorship of Mathematics and Physics at Bard College.}

\author[J. M. Jones]{James M. Jones}
\address{Department of Mathematics and Computer Science, Philipps-Universität Marburg, Hans-Meerwein-Strasse 6, D-35043 Marburg, Germany.}
\email{jonesj@mathematik.uni-marburg.de}

\author[E. Pichon-Pharabod]{Eric Pichon-Pharabod}
\address{Mathematical Institute, University of Oxford, Andrew Wiles Building, Radcliffe Observatory Quarter (550), Woodstock Road, Oxford, OX2 6GG, United Kingdom.}
\email{ eric.pichon-pharabod@maths.ox.ac.uk}
\thanks{EPP is funded by the European Union (ERC Synergy Grant MaScAmp 101167287). Views and opinions expressed are however those of the author(s) only and do not necessarily reflect those of the European Union or the European Research Council. Neither the European Union nor the granting authority can be held responsible for them.}

\author[A. Thompson]{Alan Thompson}
\address{Department of Mathematical Sciences, Loughborough University, Loughborough, Leicestershire, LE11~3TU, United Kingdom.}
\email{A.M.Thompson@lboro.ac.uk}
\thanks{For the purpose of open access, the authors have applied a Creative Commons Attribution (CC BY) licence to any Author Accepted Manuscript version arising.}

\begin{abstract}
We present some effective tools for explicit computation of N\'{e}ron-Severi lattices for generic members of families of elliptically fibred surfaces. As an application of this theory, we certify a computation of the N\'{e}ron-Severi lattice for a family of K3 surfaces associated to the tardigrade Feynman integral, originally studied by Doran, Harder, Vanhove, and Pichon-Pharabod. We also apply our tools to explicitly verify a mirror symmetry conjecture, due to Giovenzana and Thompson, in several cases, and use this to strengthen several results about the coupling group. As part of verifying this conjecture, we present a classification of Tyurin degenerations of K3 surfaces of degrees 2 and 4, and K3 surfaces polarized by the hyperbolic plane lattice.
\end{abstract}

\date{}
\maketitle

\section{Introduction}

The aims of this paper are twofold. Firstly, to develop effective tools for the certified computation of N\'{e}ron-Severi lattices in families of elliptically-fibred surfaces. Then secondly, to apply these tools to explicitly verify a mirror symmetry conjecture of L. Giovenzana and the fourth author in a number of special cases.

Our first goal arises from work of the third author. In \cite{snahlpcesp1}, he presents a certified algorithm which computes an explicit basis for the second homology $H_2(Y;\bZ)$ of an elliptically-fibred surface $\pi \colon Y \to \bP^1$ along with the bilinear form induced on it by the intersection product. This basis is described in terms of geometrically-realized cycles on $Y$. 

However, for many applications one would also like to know the N\'{e}ron-Severi lattice $\NS(Y)$, which is the sublattice generated by the algebraic cycles; this allows, for example, explicit computation of the Mordell-Weil group. \cite[Section 4.1]{snahlpcesp1} gives an algorithm to approximate this, by first computing a high precision approximation to the period mapping, then computing $\NS(Y)$ as its kernel. Unfortunately, this computation is not certified and the true N\'{e}ron-Severi lattice may be larger or smaller than the computed one.

Our main results in this direction rest on the observation that, for many purposes, one does not necessarily need to know the N\'{e}ron-Severi lattice of a specific elliptically-fibred surface; rather one is interested in the N\'{e}ron-Severi lattice of a general member of a family of surfaces. In this setting, we exploit the fact that there are certified algorithms for computing the action of monodromy on the second homology \cite{ggscs,ehpcph}. Specifically, our first main result (Proposition \ref{prop:monodromy}) shows that, if one can find suitable upper and lower bounds on the rank of the N\'{e}ron-Severi lattice $\NS(Y_t)$ of the general member $Y_t$ of a family and if the monodromy representation on $H_2(Y_t;\bZ)$ contains an irreducible subrepresentation $V$ of sufficient size, then $\NS(Y_t)$  may be computed as the orthogonal complement of $V$ in $H_2(Y_t;\bZ)$.

This approach is dependant upon finding suitable upper and lower bounds on the rank of $\NS(Y_t)$. Whilst a good lower bound can often be computed by geometric considerations (e.g. explicit computation of algebraic classes), a suitable upper bound can be more difficult to obtain. In practice, this limits the utility of this approach to families with small Mordell-Weil rank. To resolve this problem, our second main result (Proposition \ref{prop:monodromyrestriction}) extends our first result to subfamilies. This allows an inductive approach, where one starts with a large family with small Mordell-Weil rank, for which Proposition \ref{prop:monodromy} applies, then restricts to successive subfamilies. Provided that each restriction only changes the Mordell-Weil rank by a small amount, Proposition \ref{prop:monodromyrestriction} can be applied without a tight upper bound on the rank of N\'{e}ron-Severi.

As an application of these results, we study the \emph{tardigrade Feynman integral}. In \cite[Section 8]{mgtlfi}, A. Harder, P. Vanhove, and the first and third authors show that the motive of this integral is identical to that of a family of elliptically fibred K3 surfaces. In the appendix to the same paper, a numerical period computation, using the approach from \cite{snahlpcesp1}, computes an explicit approximation to the N\'{e}ron-Severi lattice of a general member of this family and estimates its rank to be $11$. Using the results of this paper, we are able to certify the result of this computation.

Our second goal comes from mirror symmetry. Along with A. Harder, the first and fourth authors have conjectured that certain degenerations of a Calabi-Yau manifold are mirror to fibration structures on its mirror Calabi-Yau manifold, and that the fibrations admit a splitting into pieces that are mirror to the components of the degenerate Calabi-Yau  \cite{mstdfcym}. This conjecture was made precise for lattice polarized K3 surfaces by L. Giovenzana and the fourth author \cite[Conjecture 5.10]{dfk3slpms}. It relates degenerations of K3 surfaces to pairs of rational surfaces glued along a smooth anticanonical curve (so-called \emph{Tyurin degenerations}), to splittings of fibrations on the mirror K3 surface. We explicitly verify this conjecture for eighteen different Tyurin degenerations polarized by three different lattices: $\langle 2 \rangle$, $\langle 4 \rangle$, and the hyperbolic plane $H$. Our main result is as follows.

\begin{theorem} \textup{(Theorem \ref{thm:mirror})}
Let $\calX \to \Delta$ be a very general Tyurin degeneration of K3 surfaces polarized by one of the lattices $\langle 2 \rangle$, $\langle 4 \rangle$, or $H$. Then a very general member $Y$ of the mirror family of lattice polarized K3 surfaces admits an elliptic fibration $\pi\colon Y \to \bP^1$ and a splitting of this fibration so that \cite[Conjecture 5.10]{dfk3slpms} holds.
\end{theorem}

This theorem represents the first explicit verification of \cite[Conjecture 5.10]{dfk3slpms} in the literature. It uses the computational tools from \cite{snahlpcesp1} and the first part of the paper to check the required lattice theoretic conditions on the elliptically fibred K3 surfaces. As part of the proof, we also classify Tyurin degenerations of $\langle 2 \rangle$-, $\langle 4 \rangle$-, and $H$-polarized K3 surfaces, extending work of Friedman \cite[Section 5]{npgttk3s} and the second and fourth authors \cite{dk3sd2,t2dk3sd4,jonesthesis}.

As an application of this Theorem, we perform an investigation into an object called the \emph{coupling group}  \cite[Definition 2.35]{dfk3slpms}. This is a somewhat mysterious object that relates the lattice polarization on a Tyurin degeneration (resp. elliptically fibred K3 surface) to the lattice polarizations on the components of the degenerate K3 surface (resp. pieces of the splitting). By studying the behaviour of the coupling group in explicit examples, we are able to strengthen several results about this group from \cite{dfk3slpms}.

\subsection{Structure of the paper}

Section \ref{sec:ellipticfibrationcomputations} concerns tools for explicitly computing the N\'{e}ron-Severi lattice of a general member of a family of elliptically fibred surfaces. We begin with Section \ref{sec:torsionMW}, which treats the simplest case where the Mordell-Weil group is known to be torsion; our main result here is Proposition \ref{prop:torsionMW}. Section \ref{sec:nontorsionMW} discusses the more difficult non-torsion case, which we approach using the monodromy representation of the family. Our main results here are Proposition \ref{prop:monodromy}, which treats the general case, and Proposition \ref{prop:monodromyrestriction}, which describes how to restrict the method to subfamilies. 

We conclude Section \ref{sec:ellipticfibrationcomputations} with some examples to illustrate these techniques. Section \ref{sec:examplecomputations} presents an explicit computation of the N\'{e}ron-Severi lattices for several families of elliptically fibred K3 surfaces polarized by lattices of high rank. The results are presented in Appendix \ref{app:elliptictables} and are used for the mirror symmetry computations in Section \ref{sec:mirror}. Finally, Section \ref{sec:tardigrade} presents the computation of the N\'{e}ron-Severi lattice for the family of K3 surfaces obtained from the tardigrade Feynman integral.

Section \ref{sec:mirror} concerns mirror symmetry. In Section \ref{sec:splittings} we begin with an overview of the relevant theory of splittings of elliptic fibrations on K3 surfaces from \cite[Section 4]{dfk3slpms}. Section \ref{sec:explicitsplitting} then applies the theory of Section \ref{sec:ellipticfibrationcomputations} and results of \cite{snahlpcesp1} to explicitly compute splittings for the fibrations studied in Section \ref{sec:examplecomputations}, along with the relevant lattice-theoretic data attached to each. The results are described in Proposition \ref{prop:splittings}.

In Section \ref{sec:tyurin} we move to the other side of the mirror and study Tyurin degenerations of K3 surfaces polarized by the lattices $\langle 2 \rangle$, $\langle 4 \rangle$, and $H$. We first summarize the relevant theory from \cite[Section 3]{dfk3slpms}, before proceeding to classify such Tyurin degenerations and study their lattice theory in Proposition \ref{prop:degenerations}. The bulk of the proof of this proposition is dedicated to the $H$-polarized case; the $\langle 2 \rangle$- and $\langle 4 \rangle$-polarized cases are essentially due to prior work of Friedman \cite[Section 5]{npgttk3s} and the second and fourth authors \cite{dk3sd2,t2dk3sd4,jonesthesis}.

In Section \ref{sec:mirrorsymmetry} we pull the results of the previous sections together to prove our main Theorem \ref{thm:mirror}, which verifies the mirror symmetry conjecture \cite[Conjecture 5.10]{dfk3slpms} in the $\langle 2 \rangle$-, $\langle 4 \rangle$-, and $H$-polarized cases. The corresponding data is summarized in Appendix \ref{app:mirrortables}. We conclude with an application to the study of the coupling group; our main results here are Corollaries \ref{cor:coupling1}, \ref{cor:coupling2}, and \ref{cor:coupling3}, which strengthen \cite[Proposition 4.27 and Lemma 5.7]{dfk3slpms}.

\medskip
\noindent\textbf{Data Access Statement.} SageMath \cite{sage} worksheets replicating all of the computations from this paper are available at \cite{DHTcode}.
\medskip

\noindent\textbf{AI and Originality Statement.} The $H$-polarized case of Theorem \ref{thm:mirror} was previously submitted to the Third Batch Benchmark of the \emph{First Proof} AI benchmarking project, which seeks to assess the performance of leading AI models against unpublished research problems. A precise statement of the problem submitted, along with our solution, can be found at \cite{firstproof}. To avoid potentially contaminating the results of the \emph{First Proof} project, AI tools were not used at any point during the preparation of this manuscript.

\section{Explicit Lattice Computations for Elliptically Fibred Surfaces} \label{sec:ellipticfibrationcomputations}

In this section we present several results that allow for certified computation of the N\'{e}ron-Severi lattice of a generic member of a family of elliptically fibred surfaces, following the approach of \cite{snahlpcesp1}.

We begin with some setup. Let $\pi\colon Y \to \bP^1$ be a relatively minimal elliptically fibred surface with fixed zero section $S_0$. We consider $H_2(Y;\bZ)$ as a lattice equipped with the usual intersection form. Let $H_{\pi} \subset H_2(Y;\bZ)$ denote the copy of the hyperbolic plane lattice generated by the classes of the zero section $S_0$ and a fibre $F$ and let $R_{\pi} \subset H_2(Y;\bZ)$ denote the sublattice generated by components of fibres that do not meet $S_0$. Following the terminology from \cite{es}, the sublattice $\Trv(Y,\pi) := H_{\pi} \oplus R_{\pi}$ is called the \emph{trivial lattice} in \cite[Definition 7]{snahlpcesp1}.

There is a well-known exact sequence of groups
\begin{equation} \label{eq:NSsequence} 0 \longrightarrow \Trv(Y,\pi) \longrightarrow \NS(Y) \longrightarrow \MW(Y,\pi) \longrightarrow 0,\end{equation}
where $\NS(Y)$ denotes the N\'{e}ron-Severi lattice of $Y$ (considered as a sublattice of $H_2(Y;\bZ)$) and $\MW(Y,\pi)$ denotes the Mordell-Weil group of the elliptic fibration $\pi\colon Y \to \bP^1$, which descends to an exact sequence
\[0 \longrightarrow R_{\pi} \longrightarrow H_{\pi}^{\perp} \longrightarrow \MW(Y,\pi) \longrightarrow 0,\]
where the orthogonal complement $H_{\pi}^{\perp}$ is taken inside $\NS(Y)$ (c.f. \cite[Section VII.2]{btes}).

The relevant results of \cite{snahlpcesp1} may be summarised as follows. Let 
\[\overline{Y} := \{y^2 = x^3 + a(u)x + b(u)\}\]
denote the Weierstrass model of $\pi \colon Y \to \bP^1$, where $u$ is an affine parameter on the base, and assume that the $a(u), b(u) \in \bQ[u]$ are polynomials with rational coefficients. Then \cite{snahlpcesp1} presents a certified algorithm which computes an explicit basis for $H_2(Y;\bZ)$, described in terms of geometrically-realized cycles on $Y$, along with the sublattices $H_{\pi}$ and $R_{\pi}$.

However, to study the geometry of $Y$, one would like to know the N\'{e}ron-Severi lattice $\NS(Y)$. \cite[Section 4.1]{snahlpcesp1} gives an algorithm to approximate this: first compute a high precision (over $100$ decimal places) approximation to the period mapping, then compute $\NS(Y)$ as its kernel. Unfortunately, this computation is not certified: in general the result may be larger or smaller than the true N\'{e}ron-Severi lattice.

In this section we present a number of results to circumvent this problem. Our approach is to realise $Y$ as a special member of an appropriate family of surfaces, then to show how the results of  \cite{snahlpcesp1} can be used to compute the correct N\'{e}ron-Severi lattice for a general member.

\subsection{Torsion Mordell-Weil in analytic families} \label{sec:torsionMW}

Our first approach works under the following assumptions. Let $\pi_t\colon {Y}_t \to \bP^1$ be a smooth family of relatively minimal elliptic surfaces with section, for $t \in \Delta^n$ (the complex polydisc), and let $\overline{Y}_t := \{y^2 = x^3 + a_t(u) x + b_t(u)\}$ be their Weierstrass models. Assume that these Weierstrass models satisfy the following assumptions.

\begin{assumption} \label{ass:analyticfamily}
\begin{enumerate}
\item $a_t(u)$ and $b_t(u)$ are polynomials in $u$ whose coefficients are holomorphic functions of $t$;
\item $a_0(u), b_0(u) \in \bQ[u]$;
\item for all $t \in \Delta^n$, the elliptic fibrations $\pi_t\colon {Y}_t \to \bP^1$ have the same singular fibre types as $\pi_0\colon Y_0 \to \bP^1$ (equivalently, the covers of $\Delta^n$ defined by the vanishing of $a_t(u)$, $b_t(u)$, and $4a_t(u)^3+27b_t(u)^2$ in $\Delta^n \times \bC$ are unramified).
\end{enumerate}
\end{assumption}

\begin{remark} The assumption of rational coefficients in Assumption \ref{ass:analyticfamily}(2) is for practical reasons, rather than mathematical ones: some of the algorithms used in \cite{snahlpcesp1} have only been implemented over $\bQ$. These algorithms should work over any number field, but there is a significant performance cost to doing so.
\end{remark}


The next result shows that, assuming the Mordell-Weil group is torsion, the N\'{e}ron-Severi lattice of a very general member $Y_t$ of the family above is just the saturation of the trivial lattice.

\begin{proposition} \label{prop:torsionMW} Suppose that Assumption \ref{ass:analyticfamily} holds and that $\MW({Y}_t,\pi_t)$ is torsion for very general $t \in \Delta^n$. Then for very general $t \in \Delta^n$, the N\'{e}ron-Severi lattice $\NS(Y_t)$ is equal to the saturation of $\Trv(Y_0,\pi_0) \cong H_{\pi} \oplus R_{\pi}$ in $H_2(Y_t;\bZ)$.
\end{proposition}

\begin{proof} Note that, under Assumption \ref{ass:analyticfamily}, $H_2(Y_t;\bZ)$ are trivial local systems over $\Delta^n$ and the trivial sublattice $\Trv(Y_t,\pi_t)$ varies continuously with $t$, so we may canonically identify $H_2(Y_t;\bZ) \cong H_2(Y_0;\bZ)$ so that $\Trv(Y_t,\pi_t)$ is identified with $\Trv(Y_0,\pi_0)$; we will do this implicitly in the remainder of the proof.

The N\'{e}ron-Severi group of $Y_t$ varies with $t$ and fits into the exact sequence \eqref{eq:NSsequence}. In the case where $\MW({Y}_t,\pi_t) \cong  \NS(Y_t)/\Trv(Y_t,\pi_t)$ is torsion, $\NS(Y_t)$ is just the saturation of $\Trv(Y_t,\pi_t)$ in $H_2({Y}_t;\bZ)$.
\end{proof}

\subsection{Non-torsion Mordell-Weil in algebraic families}\label{sec:nontorsionMW}

Our second approach allows us to consider families for which the Mordell-Weil group of the general fibre is not torsion, at the cost of some more restrictive assumptions. Let $U \subset \bA^n$ be a Zariski open set. Let $\pi_t\colon {Y}_t \to \bP^1$ be a smooth family of relatively minimal elliptic surfaces with section, for $t \in U$, and let $\overline{Y}_t := \{y^2 = x^3 + a_t(u) x + b_t(u)\}$ be their Weierstrass models. Assume that these Weierstrass models satisfy the following assumptions.

\begin{assumption}\label{ass:algebraicfamily}
\begin{enumerate}
\item $a_t(u)$ and $b_t(u)$ are polynomials in $u$ whose coefficients are rational functions of $t$;
\item There exists a point $p \in U$ such that $a_p(u), b_p(u) \in \bQ[u]$;
\item for all $t \in U$, the elliptic fibrations $\pi_t\colon {Y}_t \to \bP^1$ have the same singular fibre types as $\pi_p\colon Y_p \to \bP^1$ (equivalently, the covers of $U$ defined by the vanishing of $a_t(u)$, $b_t(u)$, and $4a_t(u)^3+27b_t(u)^2$ are unramified).
\end{enumerate}
\end{assumption}

\begin{remark}
As long as the point $p$ is sufficiently general, Assumption \ref{ass:algebraicfamily}(3) can always be arranged by shrinking the open set $U$: simply remove all loci where collisions between singular fibres occur.
\end{remark}

We further assume that bounds for the rank of the N\'{e}ron-Severi lattice of a very general member of the family $Y_t$ are known. Denote such bounds by $\rho_{\min}$ and $\rho_{\max}$, so that
\[\rho_{\min} \leq \rank(\NS(Y_t)) \leq \rho_{\max}.\]
Then the following proposition gives a criterion under which the lower bound $\rho_{\min}$ is sharp. 

\begin{proposition} \label{prop:monodromy} Suppose that Assumption \ref{ass:algebraicfamily} holds and that for very general $Y_t$:
\begin{enumerate}
\item $\rho_{\min} \leq \rank(\NS(Y_t)) \leq \rho_{\max}$;
\item $\rho_{\min} + \rho_{\max} < \rank(H_2(Y_t;\bZ)) + \rank(\Trv(Y_t,\pi_t))$; and
\item the monodromy representation of $\pi_1(U,p)$ on $\Trv(Y_p,\pi_p)^{\perp} \otimes \bQ \subset H_2(Y_p;\bQ)$ contains an irreducible subrepresentation $V$ such that 
\[\dim(V) = \rank(H_2(Y_t;\bZ)) - \rho_{\min}.\]
\end{enumerate}
Then for very general $t$ in a neighbourhood of $p \in U$, we have $V = \NS(Y_t)^{\perp} \otimes \bQ$ and $\NS(Y_t) = V^{\perp} \cap H_2(Y_t;\bZ)$. In particular, $\rank(\NS(Y_t)) = \rho_{\min}$.
\end{proposition}

\begin{remark}
In the statement of Proposition \ref{prop:monodromy}, note that the trivial lattice $\Trv(Y_p,\pi_p)$ is preserved (i.e. mapped to itself) by the action of monodromy, so $\Trv(Y_p,\pi_p)^{\perp}$ is also preserved under monodromy and therefore the monodromy representation of $\pi_1(U,p)$ on $\Trv(Y_p,\pi_p)^{\perp} \otimes \bQ$ is well-defined.
\end{remark}

\begin{proof}[Proof of Proposition \ref{prop:monodromy}]
In a simply connected neighbourhood $S$ of $p \in U$, the $H_2(Y_t;\bZ)$ define trivial local systems and the trivial sublattice $\Trv(Y_t,\pi_t)$ varies continuously with $t$, so we may identify $H_2(Y_t;\bZ) \cong H_2(Y_p;\bZ)$ so that $\Trv(Y_t,\pi_t)$ is identified with $\Trv(Y_p,\pi_p)$. We will do this implicitly in the remainder of the proof. 

Using the upper bound on the rank of $\NS(Y_t)$, we have
\[\rank(\NS(Y_t)^{\perp}) \geq \rank(H_2(Y_t;\bZ)) - \rho_{\max}.\]
Therefore
\begin{align*} \rank(\NS(Y_t)^{\perp}) + \dim(V) &\geq 2\rank(H_2(Y_t;\bZ)) - \rho_{\max} - \rho_{\min}\\
&> \rank(H_2(Y_t;\bZ)) - \rank(\Trv(Y_t,\pi_t))\\
&= \rank(\Trv(Y_t,\pi_t)^{\perp}).
\end{align*}
As $\NS(Y_t)^{\perp} \subset \Trv(Y_t,\pi_t)^{\perp}$ and $V \subset \Trv(Y_t,\pi_t)^{\perp} \otimes \bQ$, we see that $V$ must contain a nonzero element of $\NS(Y_t)^{\perp}$. Noting that $V$ is irreducible and monodromy preserves $\NS(Y_t)^{\perp}$, we must therefore have $V \subseteq \NS(Y_t)^{\perp} \otimes \bQ$. But $\dim(V) = \rank(H_2(Y_t;\bZ)) - \rho_{\min} \geq \rank(\NS(Y_t)^{\perp})$, so in fact we must have equality. The second part of the statement follows by taking orthogonal complements.
\end{proof}

To apply this proposition, typically the lower bound $\rho_{\min}$ is obtained by constructing an explicit set of divisors in $\NS(Y_t)$, which are expected to span $\NS(Y_t)$. The upper bound $\rho_{\max}$ may be obtained, for example, by computing the rank of an appropriate matrix of periods; see Section \ref{sec:tardigrade}.

Under the assumptions of Proposition \ref{prop:monodromy}, $V$ is the unique largest irreducible subrepresentation of the monodromy representation. Therefore, from a computational standpoint, to find $V$ it suffices to find a set of loops in $\pi_1(U,p)$ such that monodromy around them spans an irreducible subrepresentation of dimension~$\rank(H_2(Y_t;\bZ)) - \rho_{\min}$.

Our next result allows us to restrict this technique to subfamilies. Suppose that we have a family $Y_t$ for which Proposition \ref{prop:monodromy} holds and let $\rho := \rho_{\min} = \rank(\NS(Y_t))$. Suppose further that we have a subfamily $Y_t'$ defined over a subvariety $U' \subset U$ and a nonsingular point $p' \in U'$  such that $a_{p'}(u), b_{p'}(u) \in \bQ[u]$. Assume that the rank of the N\'{e}ron-Severi lattice of a very general member ${Y}_t'$ of the subfamily is bounded by $\rho_{\min}'$ and $\rho_{\max}'$ so that
\[\rho_{\min}' \leq \rank(\NS(Y_t')) \leq \rho'_{\max}.\]
Note that $\rho \leq \rho_{\min}'$ by semicontinuity.

Choosing a path in $U$ from $p'$ to $p$ realises the fundamental group $\pi_1(U',p')$ as a subgroup $\pi_1(U',p') \subset \pi_1(U,p)$. The monodromy representation of $\pi_1(U',p')$ on $\Trv(Y'_{p'},\pi'_{p'})^{\perp} \otimes \bQ$ (identified with $\Trv(Y_p,\pi_p)^{\perp} \otimes \bQ$ as in the proof of Proposition \ref{prop:monodromy})  is then given by restricting the monodromy representation of $\pi_1(U,p)$ to $\pi_1(U',p')$, but the subrepresentation $V$ found in Proposition \ref{prop:monodromy} will not, in general, remain irreducible.

\begin{proposition} \label{prop:monodromyrestriction} With assumptions as above, suppose that for very general $Y'_t$:
\begin{enumerate}
\item $\rho_{\min}' \leq \rank(\NS(Y_t')) \leq \rho'_{\max}$;
\item $\rho'_{\min} + \rho'_{\max} - \rho < \rank(H_2(Y_t';\bZ))$; and
\item the monodromy representation of $\pi_1(U',p')$ on $\Trv(Y'_{p'},\pi'_{p'})^{\perp} \otimes \bQ$ contains an irreducible subrepresentation $V' \subset V$ such that 
\[\dim(V') = \rank(H_2(Y_t';\bZ)) - \rho'_{\min}.\]
\end{enumerate}
Then for very general $t$ in a neighbourhood of $p' \in U'$, we have $V' = \NS(Y_t')^{\perp} \otimes \bQ$ and $\NS(Y'_t) = (V')^{\perp} \cap H_2(Y'_t;\bZ)$. In particular, $\rank(\NS(Y_t')) = \rho_{\min}'$.
\end{proposition}

\begin{proof} 
By (1) we have
\[ \rank(\NS(Y'_t)^{\perp}) \geq \rank(H_2(Y_t';\bZ)) - \rho'_{\max},\]
so using (2) and (3) we obtain
\begin{align*}\dim(V) &= \rank(H_2(Y_t';\bZ)) - \rho\\
 &<  2\rank(H_2(Y_t';\bZ)) - \rho_{\max}' - \rho_{\min}'\\
&\leq \rank(\NS(Y'_t)^{\perp}) + \dim(V').
\end{align*}
Noting that $\NS(Y_t')^{\perp} \otimes \bQ  \subset \NS(Y_t)^{\perp} \otimes \bQ = V$ and $V' \subset V$, we find that $V'$ must contain a nonzero element from $\NS(Y_t')^{\perp}$.  Noting that $V'$ is irreducible and monodromy preserves $\NS(Y_t')^{\perp}$, we must therefore have $V' \subseteq \NS(Y_t')^{\perp} \otimes \bQ$. But $\dim(V') = \rank(H_2(Y_t';\bZ)) - \rho_{\min}' \geq \rank(\NS(Y_t')^{\perp})$, so in fact we must have equality. The second part of the statement follows by taking orthogonal complements.
\end{proof}

As before, under the assumptions of Proposition \ref{prop:monodromyrestriction}, $V'$ is the unique largest irreducible subrepresentation of the action of $\pi_1(U',p')$ on $V$. Therefore, if $V$ is known, to find $V'$ it suffices to find a set of loops in $\pi_1(U',p')$ such that monodromy around them spans a subrepresentation of $V$ of dimension $\rank(H_2(Y_t';\bZ)) - \rho'_{\min}$.

The results of this section may be used to study surfaces with non-torsion Mordell-Weil group as follows. First, construct a family of Weierstrass models satisfying the conditions of Proposition \ref{prop:monodromy} which has the desired singular fibres and for which the Mordell-Weil group is expected to be small (i.e. $\MW(Y_t,\pi_t)$ is expected to have rank $\rho_{\min} - \rank(\Trv(Y_t,\pi_t))$, so one looks for a family satisfying $\rho_{\min} - \rank(\Trv(Y_t,\pi_t)) < \rank(H_2(Y_t;\bZ)) - \rho_{\max}$ in order to apply Proposition \ref{prop:monodromy}). Subsequently, restrict to subloci where the Mordell-Weil rank jumps and apply Proposition \ref{prop:monodromyrestriction}. This process can be carried out inductively if necessary to ensure that the condition $\rho'_{\min} - \rho < \rank(H_2(Y_t';\bZ)) - \rho'_{\max}$ holds at each step; note here that $\rho'_{\min} - \rho$ is the expected change in the Mordell-Weil rank.

This approach is particularly useful for K3 surfaces: it is known that for any K3 surface $Y_t$ we have $\rank(H_2(Y_t;\bZ)) = 22$ and $\rank(\NS(Y_t)) \leq 20$, so it follows that $\rank(H_2(Y_t;\bZ)) - \rho_{\max} \geq 2$. This means that the inductive approach above can be carried out for K3 surfaces without knowledge of the upper bounds $\rho_{\max}$ and $\rho'_{\max}$, as long as the expected Mordell-Weil rank of the starting family is $1$ and the expected change in Mordell-Weil rank at each step is equal to $1$.

%
%

\subsection{Example computations}\label{sec:examplecomputations}

To illustrate these methods, we explicitly compute N\'{e}ron-Severi lattices and Mordell-Weil groups for certain elliptic fibrations on lattice polarized K3 surfaces. Here we follow Alexeev and Engel's \cite{cmk3s} modification of Dolgachev's \cite{mslpk3s} original definition: Suppose that $L$ is a non-degenerate primitive sublattice of the K3 lattice $\Lambda_{\text{K3}} \cong H \oplus H \oplus H \oplus E_8 \oplus E_8$ of signature $(1,\rank({L}) -1)$. Fix a vector $h \in {L} \otimes \bR$, such that $\langle h,h \rangle_{{L}} > 0$ and $h \notin L' \otimes \bR$ for any proper sublattice $L' \subsetneq {L}$ with $\rank(L') < \rank({L})$; Alexeev and Engel call such vectors \emph{very irrational}. 

\begin{definition} \cite[Definition 2.6]{cmk3s} \label{def:Lquasipol} An \emph{${L}$-quasi\-polarized K3 surface} is a K3 surface $X$ along with a primitive lattice embedding $j\colon {L} \hookrightarrow \NS(X)$ for which $j(h)$ is nef and big.
\end{definition}

We will be particularly interested in K3 surfaces quasipolarized by the three lattices
\begin{align*}
M &:= H \oplus E_8 \oplus E_8,\\
M_1 &:= H \oplus E_8 \oplus E_8 \oplus A_1,\\
M_2 &:= H \oplus E_8 \oplus E_8 \oplus \langle -4 \rangle.
\end{align*}
Up to symmetries, there are four distinct elliptic fibrations on an $M_1$-quasipolarized K3 surface, nine on an $M_2$-quasipolarized K3 surface \cite[Remark 7.11]{mslpk3s}, and two on an $M$-quasipolarized K3 surface \cite[Section 3.3]{milpk3s}. Weierstrass models for each of these elliptic fibrations are presented in the tables in Appendix \ref{app:elliptictables}; the following examples illustrate how these were obtained.

\begin{example}\label{ex:torsionMW}
As an example to illustrate the general method, case $A_{17}$ from Table \ref{tab:M1weierstrass} was constructed as follows. First, we find the equation of a general degree $12$ curve of the form $C := \{x^3 + a(u,v) x + b(u,v) = 0\} \subset \WP(1,1,4)$, where $(u,v,x)$ have weights $(1,1,4)$ respectively, having a singularity of type $A_{17}$ on the line $\{v = 0\}$. To arrange the correct Mordell-Weil group, we further impose that there should be a curve of degree $4$ that lies tangent to $C$ to order $12$ at the singular point. The curve $C$ thus obtained is given by taking $a$ and $b$ to be the homogenizations of the polynomials $a_t(u)$ and $b_t(u)$ from Table \ref{tab:M1weierstrass}, and the tangent curve is given by
\[12t_4x - (t_1v^2 + t_2uv + t_3u^2)^2 = 0,\]
where the $t_i$ denote deformation parameters. Dehomogenizing with respect to $v$ and taking a double cover ramified along $C$ gives the Weierstrass model from Table \ref{tab:M1weierstrass}, which by construction realises a family of elliptically fibred K3 surfaces with fibres $\rmI_{18}\, 6\rmI_1$ and Mordell-Weil group $\frac{\bZ}{3\bZ}$. These K3 surfaces are quasipolarized by an index $3$ overlattice of $H \oplus A_{17}$, which is easily shown to be $M_1$. As the moduli space of $M_1$-quasipolarized K3 surfaces is $1$-dimensional and the constructed family is not isotrivial, it follows that the constructed family contains very general $M_1$-quasipolarized K3 surfaces. We may therefore apply the approach from Section \ref{sec:torsionMW} to explicitly compute their N\'{e}ron-Severi lattices.

A similar approach applies to all of the examples from Tables \ref{tab:M1weierstrass} and \ref{tab:M2weierstrass} that have torsion Mordell-Weil, with the exception of the cases that are enhancements of $D_{16}$ and $E_8E_8$ (see Example \ref{ex:Hpol}).
\end{example}

\begin{example}\label{ex:Hpol}
A modular family of $M$-quasipolarized K3 surfaces is described by Clingher and Doran \cite{milpk3s}, based on earlier work of Inose \cite{desk3sni}. This family depends on  two parameters $(t_1,t_2)$ (respectively called $(a,b)$ in \cite{milpk3s}). A very general K3 surface in this family admits two elliptic fibrations, called \emph{standard} and \emph{alternate} in \cite{milpk3s}. The standard fibration corresponds to case $E_8E_8$ in Table \ref{tab:Mweierstrass}; its Weierstrass form may be found in \cite[Section 1]{desk3sni}. The alternate fibration corresponds to case $D_{16}$ in Table \ref{tab:Mweierstrass}; its Weierstrass form may be found in \cite[Section 4.4]{milpk3s}. Given this, explicit computation of the N\'{e}ron-Severi lattices for these surfaces may be carried out using the approach from Section \ref{sec:torsionMW}.

The corresponding cases $E_8E_8A_1$ and $D_{16}A_1$ in Table \ref{tab:M1weierstrass}, along with $E_8E_8\langle -4\rangle$ and $D_{16}\langle -4 \rangle$ in Table \ref{tab:M2weierstrass} may be obtained by restricting to Noether-Lefschetz loci. Specifically, it follows from \cite[Lemma 1]{desk3sni} and \cite[Proposition 4.7]{milpk3s} that the locus where the $M$-quasipolarized fibrations $E_8E_8$ and $D_{16}$ enhance to the $M_1$-quasipolarized fibrations $E_8E_8A_1$ and $D_{16}A_1$ is given by $t_1^3 = (t_2 \pm 1)^2$, and it follows from \cite[3.2.1]{nfk3smmp} that  the locus where the $M$-quasipolarized fibrations $E_8E_8$ and $D_{16}$ enhance to the $M_2$-quasipolarized fibrations $E_8E_8\langle -4\rangle$ and $D_{16}\langle -4 \rangle$ is given by
\begin{equation} \label{eq:4polNL} \begin{array}{l} 0 = 64t_1^9 - 192t_1^6t_2^2 - 21360t_1^6 + 192t_1^3t_2^4 - 83424t_1^3t_2^2 + 1792857t_1^3  - 64t_2^6\\
\qquad -5808t_2^4 - 175692t_2^2 - 1771561. \end{array} \end{equation}
In each case, the $\bQ$ point from Tables \ref{tab:M1weierstrass} and \ref{tab:M2weierstrass} lies on the corresponding Noether-Lefschetz locus, and it follows from the analytic implicit function theorem that Assumptions \ref{ass:analyticfamily} hold in a neighbourhood of this point. In the $M_1$-quasipolarized cases (Table \ref{tab:M1weierstrass}), the Mordell-Weil group is generically torsion, so explicit computation of the N\'{e}ron-Severi lattices for these surfaces may be carried out using the approach from Section \ref{sec:torsionMW}.

In the $M_2$-quasipolarized cases (Table \ref{tab:M2weierstrass}), the Mordell-Weil group is generically rank $1$, with non-torsion part generated by a $(-4)$-class. To compute this class explicitly we use an approach from \cite{ggscs}. The K3 surface over the $\bQ$-point from Table \ref{tab:M2weierstrass} has a real structure and is thus acted upon by an antiholomorphic involution, induced by complex conjugation. The additional $(-4)$-class in the Mordell-Weil group must be a $(\pm 1)$-eigenvector for this involution and lies in the orthogonal complement of the trivial lattice, which is isometric to $H \oplus H$ and is preserved by the antiholomorphic involution. Following the approach of \cite[Section 6.2]{ggscs}, we compute the $(\pm 1)$-eigenspaces of the action of the antiholomorphic involution on $H \oplus H$. These eigenspaces contain finitely many $(-4)$-classes. Finally, we approximate the period mapping using the algorithm from \cite[Section 4.1]{snahlpcesp1}, with sufficient precision to certify that all but one of these $(-4)$-classes cannot be contained in its kernel. The remaining $(-4)$-class must be the one from the Mordell-Weil group.
\end{example}

\begin{example}\label{ex:D8D84}
The remaining case from Table \ref{tab:M2weierstrass} is $D_8D_8\langle -4 \rangle$. To construct a family of K3 surfaces of this type, we first construct an explicit family $\pi_t\colon Y_t \to \bP^1$ of (Weierstrass models of) elliptically fibred K3 surfaces with singular fibres of types $2\rmI_4^*\, 4\rmI_1$ and torsion Mordell-Weil group $\frac{\bZ}{2\bZ}$, following the same approach as Example \ref{ex:torsionMW}. Such K3 surfaces admit a second elliptic fibration, with singular fibres of types $2\rmI\rmI\rmI^*\, 2\rmI_2\, 2\rmI_1$, and one may explicitly compute the relationship between the deformation parameters in these two fibrations.

The lattice enhancement to the $D_8D_8\langle -4 \rangle$ family occurs along the locus where the two $\rmI_2$ fibres in the $2\rmI\rmI\rmI^*\, 2\rmI_2\, 2\rmI_1$ fibration collide to give an $\rmI_4$; this second fibration becomes precisely the case $E_7E_7A_3$ from Table \ref{tab:M2weierstrass}. Therefore, to construct the $D_8D_8\langle -4 \rangle$ family, it suffices to restrict the parametrization of the family $\pi_t\colon Y_t \to \bP^1$ to this locus. 

We may then apply the approach from Section \ref{sec:nontorsionMW} to explicitly compute their N\'{e}ron-Severi lattices. Here we follow the approach of \cite[Section 3]{ehpcph} and \cite{ggscs}. We first compute a set of paths generating the fundamental group, following the method from \cite[Section 3.5.3]{ehpcph}. Then we compute the Picard-Fuchs equation for $\pi_p\colon Y_p \to \bP^1$, using \cite{hfug}, and numerically integrate it with certified bounds of precision using \cite{rmpedffs} to obtain the action of monodromy on the Frobenius basis of solutions to it along these paths. From this we obtain and certify the monodromy representation on $\Trv(Y_p,\pi_p)^{\perp}$ by a change of basis. Then we apply Proposition \ref{prop:monodromy}.
\end{example}

\subsection{The tardigrade Feynman integral}\label{sec:tardigrade}

Our final example is of a somewhat different character. In \cite[Section 8]{mgtlfi}, the authors consider the \emph{tardigrade Feynman integral} and show that its motive is that of a family of elliptically fibred K3 surfaces. In the appendix to the same paper, a numerical period computation, using the approach of \cite[Section 4.1]{snahlpcesp1}, computes an explicit approximation to the N\'{e}ron-Severi lattice of a general member of this family and estimates its rank to be $11$. The following example certifies this result, using the approach from Section \ref{sec:nontorsionMW}. This illustrates the power of our methods\footnote{To quote Nigel Tufnel, ``These go to eleven''.}: K3 surfaces with Picard ranks close to the ``middle'' value of $10$ are notoriously difficult to study explicitly and, whilst this example admits several elliptic fibrations, all of the known ones have Mordell-Weil groups of rank at least $2$.

We begin with the  quartic model given in \cite[Equation (156)]{mgtlfi}. This model is a family of singular quartics with three pairs of simple nodes, given in \cite[Equation (157)]{mgtlfi}; its minimal resolution is the family of K3 surfaces we wish to study.

We first give a lower bound for the Picard rank in this family, by describing explicit algebraic classes lying on a generic member $Y_{t}$. In the notation of \cite[Equation (156)]{mgtlfi}, the two planes $\{y_0+y_1=0\}$ and $\{z_0+z_1=0\}$ each contain four nodes and the restriction of the quartic to these planes splits as a union of two conics. Moreover, the restriction of the quartic to the plane $\{y_0+y_1+z_0+z_1 = 0\}$ also splits as a union of two conics.

This gives ten classes in $\NS(Y_t)$: the hyperplane section, the six exceptional curves coming from the six nodes, and one of the conics from each of the three planes above. A further class may be constructed by considering the intersection of the quartic with the linear system of planes containing the line $\{y_0 = y_1 = 0\}$. This linear system gives rise to an elliptic fibration on $Y_t$ which has reducible singular fibres of Kodaira types: 
\begin{itemize}
\item one $\rmI_4$ fibre corresponding to the plane $\{y_0+y_1=0\}$, 
\item two $\rmI_2$ fibres corresponding to the planes through each of the two singularities with $\{z_0=z_1=0\}$, 
\item two additional $\rmI_2$ fibres. 
\end{itemize}
Define an eleventh class to be one of the components of one of the additional $\rmI_2$ fibres. Computing the rank of the intersection matrix of these $11$ classes shows that they generate a lattice of rank $11$, giving a lower bound for the Picard rank.

Next we give an upper bound on the generic Picard rank. This can be obtained from the rank of the map $H_2(Y_t;\bZ)\to \mathbb C^{n+1}$ that sends $\gamma$ to $(\int_\gamma (\nabla^{k}\omega)_t)_{0\le k\le n}$, where $n$ is an arbitrary positive integer, $\omega$ is a section of the bundle of holomorphic $2$-forms, and $\nabla$ is the action of the Gauss-Manin connection along a fixed arbitrary direction.  Indeed, by Lefschetz' Theorem on $(1,1)$ Classes, a cycle $\gamma$ is generically algebraic if and only if its image under this map is $0$ . A numerical approximation of this map for $n=11$ with certified bounds of precision is computed using \cite[Section 3.3]{snahlpcesp1}, and shows that its rank is at least $11$, and therefore the Picard rank is at most $11$.

This proves that $\rank(\NS(Y_t)) = 11$ for a generic member of the family. We can further compute the N\'{e}ron-Severi lattice explicitly inside $H_2(Y_t;\bZ)$ using Proposition~\ref{prop:monodromy}. Indeed, we may construct a second elliptic fibration $\pi_t \colon Y_t \to \bP^1$ on a generic member of the tardigrade family by considering the linear system of planes containing the line $\{y_0+y_1 = z_0 + z_1 = 0\}$. This elliptic fibration has reducible singular fibres of Kodaira types: 
\begin{itemize}
\item two $\rmI_4$'s corresponding to the planes $\{y_0+y_1=0\}$ and $\{z_0 + z_1 = 0\}$, 
\item one $\rmI_2$ corresponding to the plane $\{y_0+y_1+z_0+z_1 = 0\}$.
\end{itemize}
In particular the trivial lattice $\Trv(Y_{t},\pi_{t})$ has rank $9$. We then compute the $13$-dimensional monodromy representation on $\Trv(Y_{t},\pi_{t})^{\perp} \subset H_2(Y_t;\bQ)$ and find that it contains an irreducible subrepresentation $V$ of rank $11$. Compared to Example \ref{ex:D8D84}, the Picard-Fuchs equations that need to be integrated to compute this monodromy representation, following \cite{ehpcph,ggscs}, are too large to be computed in reasonable time using current state-of-the-art software. Instead, we computed the monodromy representation using the braid approach from \cite{vlspc,gvlses}. The braids themselves were certifiably computed using the method from \cite{cbad}, which has been implemented in Rust\footnote{\url{https://gitlab.inria.fr/aguillem/braids}}.
Finally, by Proposition~\ref{prop:monodromy}, we obtain $\NS(Y_t) = V^{\perp} \cap H_2(Y_t;\bZ)$. The resulting lattice agrees with the one computed in \cite[Appendix B]{mgtlfi}, certifying this result. A SageMath worksheet replicating this computation is available at \cite{DHTcode}.

\section{Application to mirror symmetry for fibrations and degenerations}\label{sec:mirror}

Based on work by Harder and the first and fourth authors \cite{mstdfcym}, Giovenzana and the fourth author  conjecture a mirror symmetric relationship between certain degenerations of a K3 surface and ``splittings'' of an elliptic fibration on its mirror \cite[Conjecture 5.10]{dfk3slpms}. The aim of this section is to explicitly verify the predictions of this conjecture for K3 surfaces polarized by the lattices $\langle 2 \rangle$, $\langle 4 \rangle$ and $H$ (the hyperbolic plane lattice) and their mirrors $M_1$, $M_2$, and $M$. As a consequence, we will prove sharpness for several results from \cite{dfk3slpms}.

\subsection{Splittings of elliptic fibrations}\label{sec:splittings}

The setup in this section is taken from \cite[Section 4]{dfk3slpms}. Let $Y$ be a K3 surface which admits an elliptic fibration $\pi\colon Y \to \bP^1$ with distinguished zero section $S_0$. Let $\Sigma \subset \bP^1$ denote the finite set of points over which the fibres of $\pi$ are singular.

Let $\gamma \subset \bP^1\setminus \Sigma$ be a simple loop and let $p \in \gamma$ be a point. Let $F_p = \pi^{-1}(p)$ denote the fibre over $p$ and let $F$ be its class in $H_2(Y;\bZ)$. The loop $\gamma$ divides $\bP^1$ into two pieces which are homeomorphic to open discs $\Delta_1$ and $\Delta_2$. Let $\pi_1\colon Y_1 \to \Delta_1$ and $\pi_2\colon Y_2 \to \Delta_2$ denote the induced elliptic fibrations over these discs. 

\begin{definition}\label{def:allowable} \cite[Definition 4.11]{dfk3slpms} The simple loop $\gamma \subset \bP^1 \setminus \Sigma$ is called \emph{allowable} if there exists a symplectic basis $(a,b)$ for $H_1(F_p;\bZ)$ such that:
\begin{enumerate}
\item For each $i \in \{1,2\}$, the action of anticlockwise monodromy around $\partial \Delta_i$ is given in the basis $(a,b)$ by
\[\begin{pmatrix} 1 & e(Y_i) - 12 \\ 0 & 1 \end{pmatrix},\]
where $e(Y_i)$ denotes the topological Euler number of $Y_i$.
\item The class $\overline{\tau} \in F^{\perp}/\bZ F$ obtained by parallel transporting the monodromy invariant class $a$ around $\gamma$ is primitive, where the orthogonal complement is taken inside $H_2(Y;\bZ)$ with respect to the intersection form.
\end{enumerate}
 \end{definition}
 
\begin{remark} \label{rem:allowable} We make three remarks about this definition.
\begin{enumerate}
\item When $e = 12$, the class $a \in H_1(F_p;\bZ)$ will be non-unique because there is more than one choice of monodromy invariant cycle in this case. Hence there is also more than one possibility for the class $\overline{\tau} \in H_2(Y;\bZ)$ when $e=12$.
\item Condition (2) from Definition \ref{def:allowable} is not precisely the same as the corresponding condition from \cite[Definition 4.11]{dfk3slpms}. However, the two formulations are equivalent by \cite[Lemma 2.26]{dfk3slpms} and the discussion in \cite[Section 4.3]{dfk3slpms}.
\item If condition (1) from Definition \ref{def:allowable} holds for one of $\pi_1\colon Y_1 \to \Delta_1$ or $\pi_2\colon Y_2 \to \Delta_2$, then it automatically holds for both, as the anticlockwise monodromy around $\partial \Delta_1$ is the inverse of the anticlockwise monodromy around $\partial \Delta_2$ and $e(Y_1) + e(Y_2) = e(Y) = 24$.
\end{enumerate}
\end{remark}

We henceforth assume that the loop $\gamma$ is allowable and let $(a,b)$ denote a symplectic basis of $H_1(F_p;\bZ)$ for which Definition \ref{def:allowable} holds. We now introduce a lattice polarization. Let $\check{L} \subset \Lambda_{\rmK 3}$ be a sublattice and assume that $Y$ admits an $\check{L}$-quasi\-polarization $\check{L} \hookrightarrow \NS(Y)$ in the sense of Definition \ref{def:Lquasipol}.

\begin{definition}\label{def:compatible} \cite[Definition 4.21]{dfk3slpms} The $\check{L}$-quasi\-polarization is \emph{compatible with $\pi$ and $\gamma$} if $\check{L} \subset \NS(Y)$ contains the class $F$ of a fibre of $\pi$ and the orthogonal complement $\check{L}^{\perp} \subset H^2(Y;\bZ)$ contains the class $\tau$ of a lift of $\overline{\tau}$.
\end{definition}

\begin{remark} As $\check{L}$ is non-degenerate and $F \in \check{L}$, if a lift $\tau \in \check{L}^{\perp}$ of $\overline{\tau}$ exists then it is unique, as any other lift $\tau + nF$, for some $n \neq 0$, will not be orthogonal to $\check{L}$.
\end{remark}

We will henceforth assume that the $\check{L}$-quasi\-polarization on $Y$ is compatible with $\pi$ and $\gamma$, and we will let $\tau \in \check{L}^{\perp}$ denote the class from Definition \ref{def:compatible}. Let $\Gamma$ denote the negative definite lattice $\Gamma := F^{\perp}_{\check{L}}/\bZ F$, where $F^{\perp}_{\check{L}}$ denotes the orthogonal complement of $F$ taken in $\check{L}$, and let $J = \bZ F \oplus \bZ \tau \subset H_2(Y;\bZ)$ denote the rank $2$ isotropic sublattice spanned by $F$ and $\tau$. By \cite[Lemma 4.23]{dfk3slpms} we have an embedding $\Gamma \hookrightarrow J^{\perp}/J$ and $J^{\perp}/J \cong H \oplus E_8 \oplus E_8$.

Now let $q \in \gamma$ be another point, let $F_q$ denote the fibre over $q$, and let $U := Y \setminus F_q$. Then we have an isomorphism
\begin{equation} \label{eq:homologysum} H_2(U,F_p;\bZ) \cong H_2(Y_1,F_p;\bZ) \oplus H_2(Y_2,F_p;\bZ).\end{equation}
There is a boundary map $\phi\colon H_2(U,F_p;\bZ) \to H_1(F_p;\bZ)$ defined by $\phi := \phi_1 - \phi_2$, where $\phi_i\colon H_2(Y_i,F_p;\bZ) \to H_1(F_p;\bZ)$ is the usual boundary map from the exact sequence of relative homology. 

The relative homology groups $H_2(Y_i,F_p;\bZ)$ are equipped with a (non-symmetric) bilinear form called the \emph{Seifert pairing} \cite[Definition 4.1]{dfk3slpms}. This induces a bilinear form on $\ker(\phi)$ which is symmetric, and hence imbues $\ker(\phi)$ with the structure of a lattice.

From the long exact sequence of relative homology (c.f. \cite[Lemma 4.13]{dfk3slpms}), we have an isomorphism
\[H_2(U;\bZ)/\bZ F \cong \ker(\phi) \subset H_2(U,F_p;\bZ).\]
The homology group $H_2(U;\bZ)$ is a lattice when equipped with the topological intersection form and the class $F \in H_2(U;\bZ)$ is totally isotropic, so this bilinear form descends to $H_2(U;\bZ)/\bZ F \cong \ker(\phi)$. Moreover, this bilinear form is the negative of the bilinear form induced by the Seifert pairing, and there is a surjective map of lattices $\varphi\colon H_2(U;\bZ)/\bZ F  \to F^{\perp}/\bZ F$ by \cite[Proposition 4.18]{dfk3slpms}. For each $i \in \{1,2\}$, let $r_i(a)$ denote the class in $H_2(Y_i,F_p;\bZ)$ obtained by parallel transporting $a$ around $\gamma$. 

\begin{definition}\label{def:similarclasses}
The lattice of classes in $H_2(Y_i,F_p;\bZ)$ which are \emph{similar near $F_p$} is denoted
by 
\[\SN(Y_i,F_p) := r_i(a)^{\perp}/r_i(a) = \{v \in H_2(Y_i,F_p;\bZ) \mid \phi_i(v) \text{ is a multiple of } a\}/r_i(a),\]
where the orthogonal complement $r_i(a)^{\perp}$ is taken with respect to the Seifert pairing.
\end{definition}

\begin{remark}
In \cite[Section 4.3]{dfk3slpms} the lattice $\SN(Y_i,F_p)$ is denoted by the more cumbersome notation $\NS(H_2(Y_i,F_p;\bZ))$. It follows from \cite[Theorem 4.10]{dfk3slpms} that the class $r_i(a)$ is isotropic with respect to the Seifert pairing and that the Seifert pairing is symmetric on $r_i(a)^{\perp}$, so the Seifert pairing induces a lattice structure on $\SN(Y_i,F_p)$. 
\end{remark}

Let $\Psi \subset H_2(U,F_p;\bZ)$ be the $\bZ$-submodule of rank $2$ generated by the two classes $(r_1(a),0)$ and $(0,r_2(a))$, where these expressions refer to the decomposition \eqref{eq:homologysum}. Then, by \cite[Proposition 4.20]{dfk3slpms}, $\Psi$ is contained in $\ker(\phi) \cong H_2(U;\bZ)/\bZ F$ and there is an exact sequence
\[0 \longrightarrow \bZ \longrightarrow \Psi^{\perp}/\Psi \stackrel{\varphi_Y}{\longrightarrow} J^{\perp}/J \longrightarrow 0,\]
where $\varphi_Y$ is the map induced by $\varphi$ and the orthogonal complement $\Psi^{\perp}$ is taken in $\ker(\phi) \cong H_2(U;\bZ)/\bZ F$. By \cite[Proposition 2.28]{dfk3slpms} there is an inclusion of lattices
\[\Psi^{\perp}/\Psi \subset \SN(Y_1,F_p) \oplus \SN(Y_2,F_p).\]

Define classes $r_i(b) \in H_2(Y_i,F_p;\bZ)$, for each $i \in \{1,2\}$, by parallel transporting $b$ around $\partial \Delta_i$. Note that $r_i(b)$ has nontrivial boundary $\phi_i(r_i(b)) = (12 - e(Y_i))a$. As this is a multiple of $a$, $r_i(b)$ descends to a class in $\SN(Y_i,F_p)$; denote this class by $[r_i(b)]$. We now summarise several properties of $[r_i(b)]$, most of which are originally due to \cite{dfk3slpms}. 

\begin{lemma}\label{lem:r_i(b) properties}
The classes $[r_i(b)]$ have the following properties:
\begin{enumerate}[(i)]
    \item The class $([r_1(b)], -[r_2(b)]) \in  \SN(Y_1,F_p) \oplus \SN(Y_2,F_p)$ is contained in the kernel $\ker(\varphi_Y) \subset \Psi^{\perp}/\Psi$ .
\item For each $i \in \{1,2\}$ we have 
\begin{equation}\label{eq:ribperp} [r_i(b)]^{\perp} = (\Psi^{\perp} \cap H_2(Y_i,F_p;\bZ))/\bZ r_i(a),\end{equation}
where the orthogonal complement $[r_i(b)]^{\perp}$ is taken in $\SN(Y_i,F_p)$. Moreover, we have $\Psi^{\perp} \cap H_2(Y_i,F_p;\bZ) = \ker(\phi_i)$, so
\[[r_i(b)]^{\perp} = \ker(\phi_i)/\bZ r_i(a).\]
\item The map $\varphi_Y$ is injective on $[r_i(b)]^{\perp}$, for each $i \in \{1,2\}$.
\end{enumerate}
\end{lemma}
\begin{proof}
These results are essentially restatements of general results about pseudolattices from \cite[Section 2.4]{dfk3slpms}, using the identifications of the various objects in \cite[Section 4.3]{dfk3slpms}. Specifically, (i) follows from  \cite[Proposition 2.28 and Lemma 2.29]{dfk3slpms} and (iii) follows from \cite[Lemma 2.33]{dfk3slpms}.

The only part requiring further clarification is (ii). From \cite[Lemma 2.32]{dfk3slpms} we obtain
\[  (\Psi^{\perp}/\Psi) \cap \SN(Y_i,F_p) = [r_i(b)]^{\perp}.\]
Combining this with the fact that
\[\Psi^{\perp} \subset \ker(\phi) \subset H_2(Y_1,F_p;\bZ) \oplus H_2(Y_2,F_p;\bZ)\]
and with the definition of $\SN(Y_i,F_p)$ gives Equation \eqref{eq:ribperp}.

It thus remains to show that $\Psi^{\perp} \cap H_2(Y_i,F_p;\bZ) = \ker(\phi_i)$. The inclusion $\Psi^{\perp} \cap H_2(Y_i,F_p;\bZ) \subseteq \ker(\phi_i)$ is immediate: $\Psi^{\perp}$ is contained in $\ker(\phi)$ by definition and the restriction of $\phi$ to $H_2(Y_i,F_p;\bZ)$ is $(-1)^{i-1}\phi_i$. In the other direction, $\ker(\phi_i) \subseteq \ker(\phi)$ so it suffices to show that $\ker(\phi_i)$ is orthogonal to $\Psi$. This follows from the adjoint property of $\phi_i$ and $r_i$: for any $v \in \ker(\phi_i)$ and $w = w_1r_1(a) + w_2r_2(a) \in \Psi$ we have
\[\langle v,w\rangle_{\ker(\phi)} = \langle v, w_ir_i(a) \rangle_{H_2(Y_i,F_p;\bZ)} = \langle \phi_i(v), w_ia \rangle_{H_1(F;\bZ)} = \langle 0, w_ia \rangle_{H_1(F;\bZ)} = 0,\]
where the subscripts denote the $\bZ$-modules in which the pairings are evaluated.
\end{proof}

We thus obtain three sublattices of $J^{\perp}/J \cong H \oplus E_8 \oplus E_8$ associated to an $\check{L}$-quasipolarized K3 surface $\pi\colon Y \to \bP^1$ with an allowable loop $\gamma \subset \bP^1$, where $\check{L}$ is compatible with both $\pi$ and $\gamma$. Namely, we obtain: 

\begin{itemize}
\item The lattice $\Gamma := F^{\perp}_{\check{L}}/\bZ F$;
\item The lattices $\varphi_Y([r_i(b)]^{\perp})$, for $i \in \{1,2\}$.
\end{itemize}

By \cite[Proposition 5.3]{dfk3slpms}, it is significantly simpler to verify \cite[Conjecture 5.10]{dfk3slpms} under an additional assumption.

\begin{definition} \label{def:1admissible} \cite[Section 5]{mslpk3s} Let $L$ be a nondegenerate lattice and $e \in L$ be a primitive isotropic element. Then $e$ is called \emph{$1$-admissible in $L$} if there exists an isotropic class $g \in L$ with $\langle e,g\rangle_{L} = 1$.
\end{definition}

\begin{definition} \label{def:fibrationdoublyadmissible} \cite[Definition 4.29]{dfk3slpms} An $\check{L}$-quasi\-polarization on $Y$ that is compatible with $\pi$ and $\gamma$ is called \emph{doubly admissible} if $\tau$ is $1$-admissible in $\check{L}^{\perp}$ and $F$ is $1$-admissible in $\check{L}$, in the sense of Definition \ref{def:1admissible}.
\end{definition}

\begin{remark}This assumption is quite natural, for the following reasons:
\begin{itemize}
\item $1$-admissibility of $F$ in $\check{L}$ is essentially equivalent to the condition that the class of the section $S_0$ is contained in the polarizing lattice $\check{L} \subset \NS(Y)$ \cite[Proposition 4.30]{dfk3slpms}.
\item $1$-admissibility of $\tau$ in $\check{L}^{\perp}$ is required for the mirror correspondence to work nicely; see \cite[Section 6]{mslpk3s} and \cite[Remark 5.12]{dfk3slpms}.
\item $1$-admissibility of $\tau$ in $\check{L}^{\perp}$ also ensures that the sublattice $\Gamma \subset J^{\perp}/J$ is primitive \cite[Proposition 4.28]{dfk3slpms}.
\end{itemize}
\end{remark}

\subsection{Explicit computation of splitting data}\label{sec:explicitsplitting}

By \cite[Proposition 5.3]{dfk3slpms}, in order to verify \cite[Conjecture 5.10]{dfk3slpms} we need to find a simple loop $\gamma$, verify that it is allowable and compatible with the lattice polarization, check that the polarization is doubly admissible, and explicitly compute the three sublattices $\Gamma$, $\varphi_Y([r_1(b)]^{\perp})$ and $\varphi_Y([r_2(b)]^{\perp})$ of $J^{\perp}/J \cong H \oplus E_8 \oplus E_8$. We will do this for each of the elliptic fibrations on K3 surfaces quasipolarized by the lattices $M_1$, $M_2$, and $M$, as listed in Tables \ref{tab:Mweierstrass}, \ref{tab:M1weierstrass}, and \ref{tab:M2weierstrass}, via the following method.
\begin{method}\label{met:computesplittingdata}~
\begin{enumerate}
\item Let $\pi\colon Y \to \bP^1$ be an elliptically fibred K3 surface in Weierstrass form with rational coefficients, which is the central fibre of a family satisfying Assumptions \ref{ass:analyticfamily}.  The algorithm from \cite[Definition 7]{snahlpcesp1} allows us to compute an explicit basis for the lattice $H_2(Y;\bZ)$, along with the sublattices $H_{\pi}$ and $R_{\pi}$ with respect to this basis. Assume that we may find an explicit sublattice $\check{L} \subset H_2(Y;\bZ)$ which is equal to the N\'{e}ron-Severi lattice of a very general member of the family containing $Y$ (for example, by using the methods from Section \ref{sec:ellipticfibrationcomputations}).
\item Use linear algebra to compute an explicit basis for $H^{\perp}_{\pi} \subset H_2(Y;\bZ)$ and express the bilinear form in this basis. As there is a unique embedding $H \subset \Lambda_{\rmK 3}$ up to isometry \cite[Proposition 1.6.1]{isbfa}, $H^{\perp}_{\pi}$ is isometric to $H \oplus H \oplus E_8 \oplus E_8$.  The composition $H_{\pi}^{\perp} \hookrightarrow F^{\perp} \to F^{\perp}/\bZ F$ is an isometry of lattices, so we may identify $H_{\pi}^\perp \cong F^{\perp}/\bZ F$.
\item Part of the algorithm from \cite{snahlpcesp1} computes the action of monodromy around simple loops around each singular fibre in $\pi\colon Y \to \bP^1$, with respect to a base point $p \in \bP^1 \setminus \Sigma$. Using this, search for classes of simple loops $\gamma$ in $\pi_1(\bP^1 \setminus \Sigma,p)$ for which Definition \ref{def:allowable}(1) holds.
\item For each such loop $\gamma$, compute the class $\overline{\tau} \in F^{\perp}/\bZ F$ by parallel transporting the monodromy invariant class $a$ around $\gamma$ and check that it is primitive. This verifies Definition \ref{def:allowable}(2) and proves that the loop $\gamma$ is allowable.
\item Compute the lift $\tau$ of $\overline{\tau}$ to $H_{\pi}^{\perp}$ in terms of the basis from (2).  Check that $\tau$ is  contained in the orthogonal complement $\check{L}^{\perp}$. This verifies that the $\check{L}$-quasipolarization is compatible with $\pi$ and $\gamma$.
\item Check $\tau$ is $1$-admissible in $\check{L}^{\perp}$ by finding a class $g \in \check{L}^{\perp}$ such that $\tau$ and $g$ generate a copy of the hyperbolic plane lattice. Denote the lattice generated by $\tau$ and $g$ by $H_{\tau}$; note that $H_{\tau} \subset H_{\pi}^{\perp}$.
\item Part of the algorithm from \cite{snahlpcesp1} computes an explicit basis for the relative homology $H_2(U,F_p;\bZ)$. This relative homology group splits as a direct sum of relative homologies in a neighbourhood of each singular fibre, allowing us to explicitly split $H_2(U,F_p;\bZ)$ as a direct sum of relative homology groups $H_2(Y_i,F_p;\bZ)$, realising the isomorphism from Equation \eqref{eq:homologysum}. 
\item Explicitly compute the boundary maps $\phi_i\colon H_2(Y_i,F_p;\bZ) \to H_1(F_p;\bZ)$ and hence compute the kernels $\ker(\phi_i)$, for each $i \in \{1,2\}$, as sublattices of $H_2(U,F_p;\bZ)$.
\item Composing the map $\varphi$ with the isomorphism $H_{\pi}^\perp \cong F^{\perp}/\bZ F$ gives a surjective map
\[\overline{\varphi}\colon \ker(\phi) \longrightarrow H_{\pi}^{\perp}.\]
Compute the images $\overline{\varphi}(\ker(\phi_i)) \subset H_{\pi}^{\perp}$, for each $i \in \{1,2\}$.
\item Using linear algebra, find an explicit basis for $J_{\tau} := (H_{\pi} \oplus H_{\tau})^{\perp} \subset H_2(Y;\bZ)$ and express the bilinear form in this basis. As there is a unique embedding $H \oplus H \subset \Lambda_{\rmK 3}$ up to isometry \cite[Proposition 1.6.1]{isbfa}, $J_{\tau}$ is isometric to $H \oplus E_8 \oplus E_8$. \cite[Lemma 5.4]{mslpk3s} implies that we may decompose
\[H_2(Y;\bZ) = H_{\pi} \oplus H_{\tau} \oplus J_{\tau},\]
from which it follows that there is an isometry $J_{\tau} \cong J^{\perp}/J$.
\item Equation \eqref{eq:ribperp} implies that $\overline{\varphi}(\ker(\phi_i))$ is orthogonal to $\tau$. Compute the sublattice $g_i^{\perp}$ of $\overline{\varphi}(\ker(\phi_i))$ consisting of elements orthogonal to $g$. Then $g_i^{\perp}$ is contained in $J_{\tau}$ and it follows from Lemma \ref{lem:r_i(b) properties} that $g_i^{\perp}$ is identified with $\varphi_Y([r_i(b)]^{\perp})$ under the isometry $J_{\tau} \cong J^{\perp}/J$.
\item Finally, compute $\Gamma := \check{L} \cap J_{\tau}$. By construction, we will have $\Gamma \cong F^{\perp}_{\check{L}}/\bZ F$.
\end{enumerate}
\end{method}

\begin{remark}
Recall from Remark \ref{rem:allowable}(1) that when $e = 12$ the class $a \in H_1(F_p;\bZ)$ will be non-unique, as in this case there is more than one choice of monodromy invariant cycle, and thus there will also be more than one possibility for the class $\tau$. If this occurs, steps (5) and (6) provide constraints on the possible $a$'s that can be chosen.
\end{remark}

Note that the computations above take place on the K3 surface $Y$, but the lattice $\check{L} \subset H_2(Y;\bZ)$ is only known to be equal to the N\'{e}ron-Severi lattice for a very general member of the corresponding family of K3 surfaces. The following lemma shows that the computation above extends to a neighbourhood of $Y$, and hence gives the correct lattices $\Gamma$, $\varphi_Y([r_1(b)]^{\perp})$ and $\varphi_Y([r_2(b)]^{\perp})$ for a very general member of the family.

\begin{lemma} \label{lem:deformation} Assume that there exists a family of elliptically fibred K3 surfaces $\pi_t\colon Y_t \to \bP^1$ satisfying Assumption \ref{ass:analyticfamily}, such that $\pi\colon Y \to \bP^1$ is the member of the family corresponding to $t = 0$. Assume further that the sublattice $\check{L} \subset H_2(Y;\bZ)$ is equal to the N\'{e}ron-Severi lattice of a very general member of this family. Then there exists an open (analytic) neighbourhood $U \subset \Delta^n$ of $t=0$ such that Method \ref{met:computesplittingdata} computes $\Gamma$, $\varphi_{Y_t}([r_1(b)]^{\perp})$ and $\varphi_{Y_t}([r_2(b)]^{\perp})$ for all ${Y}_t$ with $t \in U$.
\end{lemma}

\begin{remark}
Note that Assumption \ref{ass:analyticfamily} is strictly weaker than Assumption \ref{ass:algebraicfamily}, so this result also holds when we have computed the N\'{e}ron-Severi lattice using the approach from Subsection \ref{sec:nontorsionMW}.
\end{remark}

\begin{proof} Recall that $\Sigma \subset \bP^1$ denotes the set of points over which the fibres of $\pi\colon Y \to \bP^1$ are singular. Let $\delta > 0$ and let 
\[\Sigma_{\delta} := \{x \in \bP^1 : | x - s| < \delta \text{ for some } s \in \Sigma\}\]
be a $\delta$-neighbourhood of $\Sigma$. Choose $\delta > 0$ sufficiently small such that $\Sigma_{\delta}$ does not intersect the loop $\gamma$ and each connected component of $\Sigma_{\delta}$ contains precisely one point from $\Sigma$. 

The singular fibre type of ${Y}_t$ is independent of $t$ and the singular fibres occur at roots of the discriminant $4a_t(u)^3 + 27b_t(u)^2$, which is a polynomial in $u$ varying holomorphically with $t$. Therefore we may find a neighbourhood $U$ of $t=0$ such that for all $t \in U$, the singular fibres of ${Y}_t$ occur over points in $\Sigma_{\delta}$. Then the loop $\gamma$ found in step (3) of Method \ref{met:computesplittingdata} satisfies the required monodromy properties for all $t \in U$ and the groups $H_2(({Y}_t)_i,F_p;\bZ)$ are independent of $t$. The remaining computations find continuously varying classes or sublattices in the trivial linear systems $H_2(({Y}_t)_i,F_p;\bZ)$ and $H_2({Y}_t;\bZ)$, and are therefore independent of $t$.
\end{proof}

By applying Method \ref{met:computesplittingdata} to the elliptic fibrations from Tables \ref{tab:Mweierstrass}, \ref{tab:M1weierstrass}, and \ref{tab:M2weierstrass}, we obtain the following result. SageMath worksheets replicating these computations are available at \cite{DHTcode}.

\begin{proposition}\label{prop:splittings}
Up to symmetry, the elliptic fibrations from Tables \ref{tab:Mweierstrass}, \ref{tab:M1weierstrass}, and \ref{tab:M2weierstrass} admit allowable splitting loops, which are compatible with the quasipolarizations as listed in Table \ref{tab:splittingdata}. All of the quasipolarizations are doubly admissible.
\end{proposition}

\begin{remark} Step (3) in Method \ref{met:computesplittingdata} was carried out by searching for loops in $\bP^1 \setminus \Sigma$ with the correct monodromy. This process was not exhaustive, so we cannot guarantee that we have found \emph{all} possible allowable splittings of the elliptic fibrations from Tables \ref{tab:Mweierstrass}, \ref{tab:M1weierstrass}, and \ref{tab:M2weierstrass}, although we expect this to be the case. \end{remark}

\subsection{Tyurin degenerations of K3 surfaces}\label{sec:tyurin}

The setup in this section is taken from \cite[Section 3]{dfk3slpms}. 

\begin{definition} \label{def:tyurin}
An \emph{$L$-polarized Tyurin degeneration of K3 surfaces} is a semistable degeneration of K3 surfaces $\pi \colon \calX \to \Delta$, where $\Delta$ denotes the open complex unit disc, such that
\begin{enumerate}
\item the central fibre $X_0$ is the union $V_1 \cup_C V_2$ of two rational surfaces $V_1$ and $V_2$ glued along a smooth curve $C$ that is anticanonical in each (i.e. $C \in |-K_{V_1}|$ and $C \in |-K_{V_2}|$),
\item (monodromy invariance) over the punctured disc $\Delta^*$ there is a trivial local subsystem $\calL$ of $R^2\pi_*\bZ$ with fibre $L$, so that for each $p \in \Delta^*$ the fibre $\calL_p \subset H^2(X_p;\bZ)$ defines an $L$-quasi\-polarization of the fibre $X_p$ over $p$,
\item (positivity) there exists an effective, flat, and nef Cartier divisor $A$ on $\calX$  such that the restriction $A_p$ of $A$ to any fibre $X_p$ over $p \in \Delta^*$ is nef and big and satisfies $[A_p] \in \calL_p$. Moreover, the restriction of $A$ to $V_i$ satisfies $(A|_{V_i})^2 > 0$ for each $i \in \{1,2\}$.
\end{enumerate}
\end{definition}

We further assume that our Tyurin degeneration is very general, in the following sense. Equipping $\calX \to \Delta$ with a choice of mixed marking \cite[Definition 4.15]{cmk3s}, by the results of \cite[Section 4]{cmk3s} there is a period map from $\Delta$ to an appropriate partial compactification of the moduli space of partially marked $L$-quasipolarized K3 surfaces. A Tyurin degeneration is \emph{very general} if the image of $0 \in \Delta$ under this period map does not lie in the closure of any Noether-Lefschetz locus.

We denote a general fibre in an $L$-polarized Tyurin degeneration by $X$. By \cite[Lemma 3.5]{npgttk3s}, the weight filtration on $H^2(X)$ given by the limiting mixed Hodge structure is defined over $\bZ$ and takes the form
\[0 \subset I \subset I^{\perp} \subset H^2(X),\]
where $I$ is a primitive rank $2$ isotropic sublattice of $H^2(X;\bZ)$. It is a consequence of \cite[Lemma 4.1.2]{cmsak3s} that there is a unique such $I$ up to isometries of $H^2(X;\bZ)$, so the second weight graded piece $I^{\perp}/I$ is always isometric to $H \oplus E_8 \oplus E_8$.

By \cite[Lemma 3.22]{dfk3slpms}, the polarizing lattice $L$ is orthogonal to $I$ and we have an embedding of $L$ into $I^{\perp}/I$. We may therefore treat $L$ as a sublattice of $I^{\perp}/I$.

\begin{definition} \cite[Definition 3.31]{dfk3slpms} An $L$-polarized Tyurin degeneration is called \emph{doubly admissible} if there exist primitive generators $e_1$ and $e_2$ for $I$ and isotropic classes $g_1$ and $g_2$ in $L^{\perp}$ such that $\langle e_1,g_1 \rangle = \langle e_2,g_2\rangle = 1$ and $\langle e_1,g_2 \rangle = \langle e_2,g_1 \rangle = 0$.
\end{definition}

\begin{remark} \label{rem:doublyadmissible} By \cite[Lemma 3.32]{dfk3slpms}, the doubly admissible condition is equivalent to the primitive embedding $I \hookrightarrow L^{\perp}$ factoring through primitive embeddings $I \hookrightarrow H \oplus H \hookrightarrow L^{\perp}$.
\end{remark}

One consequence of this definition is that if the $L$-polarized Tyurin degeneration is doubly admissible, then the sublattice $L \subset I^{\perp}/I$ is primitive.

Next, consider the class $\xi \in H^2(V_1;\bZ) \oplus H^2(V_2;\bZ)$ defined by $\xi := (C,-C) = (-K_{V_1},K_{V_2})$. As $K_{V_1}^2 = -K_{V_2}^2$ for a Tyurin degeneration, $\xi$ is isotropic. By \cite[Lemma 3.5]{npgttk3s}, the second weight graded piece of the Clemens-Schmid sequence induces an exact sequence of lattices
\[0 \longrightarrow \bZ \xi \longrightarrow \xi^{\perp} \stackrel{\varphi_X}{\longrightarrow} I^{\perp}/I \longrightarrow 0.\]
It follows from \cite[Lemma 2.32 and Proposition 3.17]{dfk3slpms} that $\xi^{\perp} \cap H^2(V_i;\bZ) = K_{V_i}^{\perp}$, and by \cite[Lemma 2.33]{dfk3slpms} the map $\varphi_X$ is injective on $K_{V_1}^{\perp}$ and $K_{V_2}^{\perp}$. We therefore obtain three sublattices of $I^{\perp}/I \cong H \oplus E_8 \oplus E_8$:
\begin{itemize}
\item The lattice $L$;
\item The lattices $\varphi_X(K_{V_i}^{\perp})$, for $i \in \{1,2\}$.
\end{itemize}
The aim of this section is to compute these three lattices for Tyurin degenerations polarized by the lattices $\langle 2 \rangle$, $\langle 4 \rangle$, and $H$.

We now introduce notation to describe these sublattices explicitly. By \cite[Lemma 1.14]{tiidk3s}, the components $V_1$ and $V_2$ of $X_0$ are either $\bP^1 \times \bP^1$, the Hirzebruch surface $\bF_2$, or a blowup $\mathrm{Bl}_n(\bP^2)$ of $\bP^2$ in $0 \leq n \leq 18$ points. We begin by describing explicit bases $\{f_i\}$ for $H^2(V_i;\bZ)$ in each case.
\medskip

\textbf{Case $\bP^1 \times \bP^1$ or $\bF_2$.} When $V_i \cong \bP^1 \times \bP^1$, let $a,b$ denote the standard basis of $H^2(\bP^1 \times \bP^1;\bZ)$ arising from the two rulings. In this basis we have 
\[K_{V_i} = -2a-2b\] 
and we define classes $f_1 = a-b$ and $f_2 = b$.

In the case where $V_i \cong \bF_2$, let $s,f$ denote the standard basis of $H^2(\bF_2;\bZ)$ given by the class $s$ of the $(-2)$-section and a fibre $f$ of the ruling. In this basis we have 
\[K_{V_i} = -2s-4f\] 
and we define classes $f_1 = s$ and $f_2 = f$.

\begin{figure}
\begin{tikzpicture}[scale=1,thick]
\draw (3,0)--(3,1);
\draw (1,1)--(4,1);
\draw[dashed] (4,1)--(5.5,1);
\draw (5.5,1)--(6.5,1);

\node[below] at (1,0.9) {$f_1$};
\node[below] at (2,0.9) {$f_2$};
\node[below right] at (2.9,0.9) {$f_3$};
\node[below] at (4,0.9) {$f_4$};
\node[below] at (5.5,0.9) {$f_{n-1}$};
\node[below] at (6.5,0.9) {$f_{n}$};

\node[right] at (3,0) {$f_0$};

\node[above] at (1,1.1) {$-2$};
\node[above] at (2,1.1) {$-2$};
\node[above] at (3,1.1) {$-2$};
\node[above] at (4,1.1) {$-2$};
\node[above] at (5.5,1.1) {$-2$};
\node[above] at (6.5,1.1) {$-1$};

\node[left] at (2.9,0) {$-2$};

\filldraw[black]
(3,0) circle (3pt)
(2,1) circle (3pt)
(1,1) circle (3pt)
(3,1) circle (3pt)
(4,1) circle (3pt)
(5.5,1) circle (3pt)
(6.5,1) circle (3pt)
;

\node at (3.75,-0.7) {$\mathrm{Bl}_n(\bP^2)$ with $n \geq 3$};

\draw (7,3)--(7,4);
\draw (6,4)--(7,4);

\node[below] at (6,3.9) {$f_1$};
\node[below right] at (6.9,3.9) {$f_2$};
\node[right] at (7,3) {$f_0$};

\node[above] at (6,4.1) {$-2$};
\node[above] at (7,4.1) {$-1$};
\node[left] at (6.9,3) {$-1$};

\filldraw[black]
(7,3) circle (3pt)
(6,4) circle (3pt)
(7,4) circle (3pt)
;

\node at (6.5,2.275) {$\mathrm{Bl}_2(\bP^2)$};

\draw (3.5,4)--(4.5,4);

\node[below] at (3.5,3.9) {$f_1$};
\node[below] at (4.5,3.9) {$f_2$};

\node[above] at (3.5,4.1) {$-2$};
\node[above] at (4.5,4.1) {$0$};

\filldraw[black]
(3.5,4) circle (3pt)
(4.5,4) circle (3pt)
;

\node at (4,2.3) {$\bP^1 \times \bP^1$ or $\bF_2$};

\draw (2,3)--(2,4);

\node[right] at (2,4) {$f_1$};
\node[right] at (2,3) {$f_0$};

\node[left] at (1.9,4) {$-1$};
\node[left] at (1.9,3) {$0$};

\filldraw[black]
(2,4) circle (3pt)
(2,3) circle (3pt)
;

\node at (2,2.25) {$\bF_1$};

\node[right] at (0.5,3) {$f_0$};

\node[left] at (0.4,3) {$1$};

\filldraw[black]
(00.5,3) circle (3pt)
;

\node at (0.5,2.3) {$\bP^2$};

\draw (0,2)--(7.5,2);
\draw (1.25,2)--(1.25,4.5);
\draw (2.75,2)--(2.75,4.5);
\draw (5.25,2)--(5.25,4.5);

\end{tikzpicture}
\caption{Dual graphs of the classes $f_i$; the numbers denote self-intersections.}
\label{fig:dualgraphs}
\end{figure}
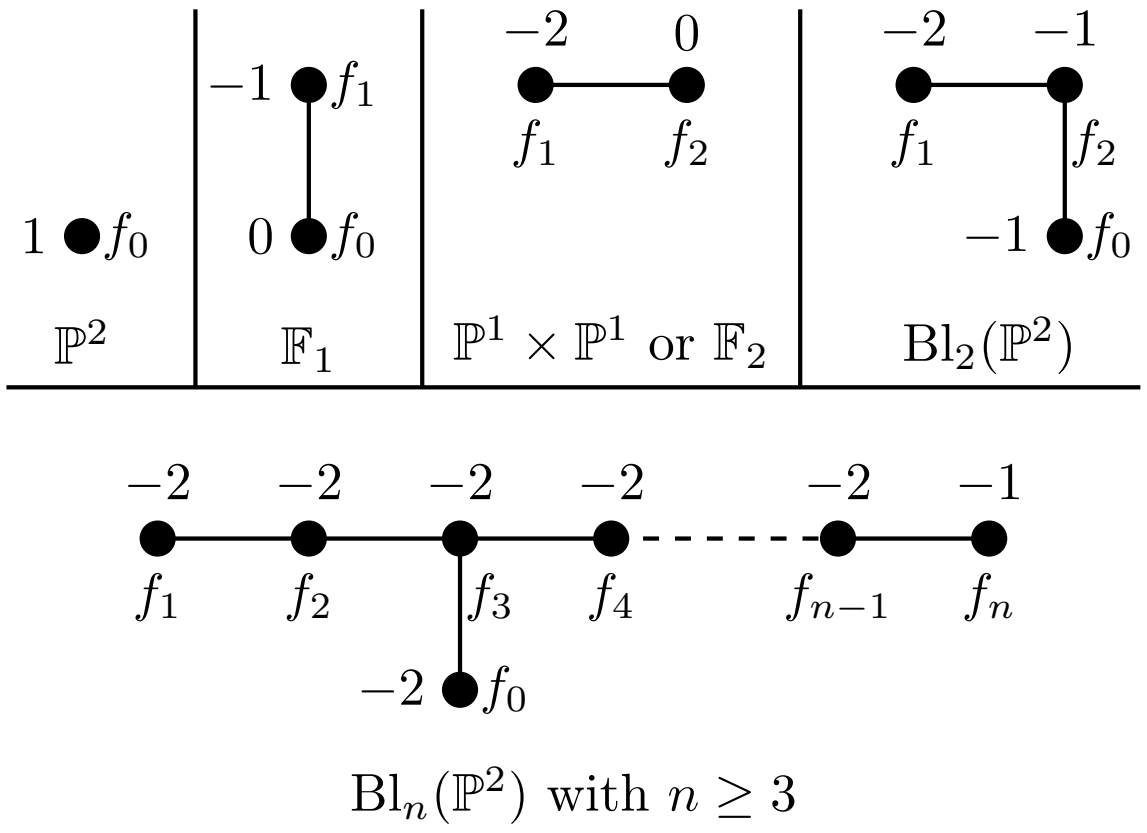

In both cases $\{f_1,f_2\}$ is a $\bZ$-basis for $H^2(V_i;\bZ)$; its dual intersection graph is shown in Figure \ref{fig:dualgraphs}. In terms of this basis we have $K_{V_i} = -2f_1-4f_2$ and $K_{V_i}^{\perp} = \langle f_1 \rangle$.
\medskip

\textbf{Case $\mathrm{Bl}_n(\bP^2)$.} When $V_i \cong \mathrm{Bl}_n(\bP^2)$,  let $l,e_1,\ldots,e_{n} \in H^2(V_i;\bZ)$ denote the standard basis arising from the total transform $l$ of a line in $\bP^2$ and the exceptional $(-1)$-class $e_i$ arising from the $i$-th blow up. In this basis we have 
\[K_{V_i} = -3l + \sum_{i=1}^{n} e_i\] 
and define classes $f_j$ as follows:
\begin{itemize}
\item If $n = 0$, define $f_0 = l$.
\item If $n = 1$, define $f_0 = l - e_1$ and $f_1 = e_1$.
\item If $n = 2$, define $f_0 = l - e_1-e_2$, $f_1 = e_1 - e_2$, and $f_2 = e_2$.
\item If $n \geq 3$, define $f_0 = l-e_1-e_2-e_3$, $f_i = e_i - e_{i+1}$ for $1 \leq i \leq n-1$, and $f_{n} = e_{n}$.
\end{itemize} 
In each case $\{f_0,\ldots,f_{n}\}$ is a $\bZ$-basis for $H^2(V_i;\bZ)$; its dual intersection graph is shown in Figure \ref{fig:dualgraphs}. The orthogonal complements $K_{V_i}^{\perp}$ are easily computed:
\begin{itemize}
\item If $V_i \cong \bP^2$, then $K_{V_i}^{\perp} = \{0\}$.
\item If $V_i \cong \mathrm{Bl}_1(\bP^2) = \bF_1$, then $K_{V_i}^{\perp} = \langle 2f_1 - f_0 \rangle$.
\item If $V_i \cong \mathrm{Bl}_2(\bP^2)$, then $K_{V_i}^{\perp} = \langle f_1, f_2 - f_0 \rangle$.
\item If $V_i \cong \mathrm{Bl}_n(\bP^2)$ for $n \geq 3$, then $K_{V_i}^{\perp} = \langle f_0,f_1,\ldots,f_{n-1} \rangle$.
\end{itemize}
\medskip

\begin{figure}
\begin{tikzpicture}[scale=0.7,thick]
\draw (2,0)--(2,1);
\draw (0,1)--(16,1);
\draw(14,0)--(14,1);

\node[above] at (0,1) {$r_1$};
\node[above] at (1,1) {$r_2$};
\node[above] at (2,1) {$r_3$};
\node[above] at (3,1) {$r_4$};
\node[above] at (4,1) {$r_5$};
\node[above] at (5,1) {$r_6$};
\node[above] at (6,1) {$r_7$};
\node[above] at (7,1) {$r_8$};
\node[above] at (8,1) {$r_9$};
\node[above] at (9,1) {$r_{10}$};
\node[above] at (10,1) {$r_{11}$};
\node[above] at (11,1) {$r_{12}$};
\node[above] at (12,1) {$r_{13}$};
\node[above] at (13,1) {$r_{14}$};
\node[above] at (14,1) {$r_{15}$};
\node[above] at (15,1) {$r_{16}$};
\node[above] at (16,1) {$r_{17}$};

\node[left] at (2,0) {$r_0$};
\node[right] at (14,0) {$r_{18}$};

\filldraw[black]
(2,0) circle (3pt)
(14,0) circle (3pt)
(1,1) circle (3pt)
(0,1) circle (3pt)
(2,1) circle (3pt)
(3,1) circle (3pt)
(4,1) circle (3pt)
(5,1) circle (3pt)
(6,1) circle (3pt)
(7,1) circle (3pt)
(8,1) circle (3pt)
(9,1) circle (3pt)
(10,1) circle (3pt)
(11,1) circle (3pt)
(12,1) circle (3pt)
(13,1) circle (3pt)
(14,1) circle (3pt)
(15,1) circle (3pt)
(16,1) circle (3pt)
;

\end{tikzpicture}
\caption{Dual graph of the $(-2)$-classes $r_i$.}
\label{fig:dualgraph}
\end{figure}
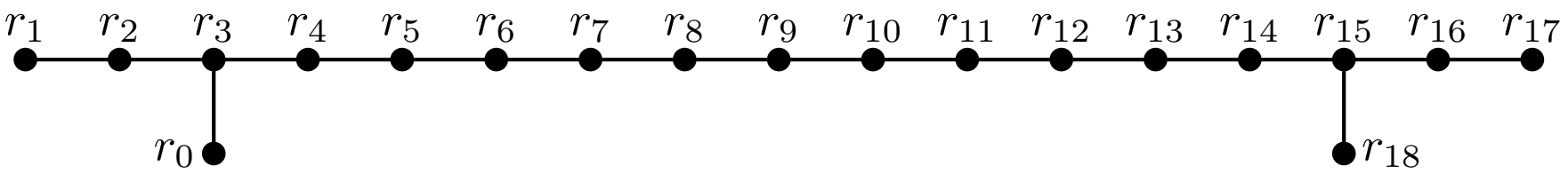

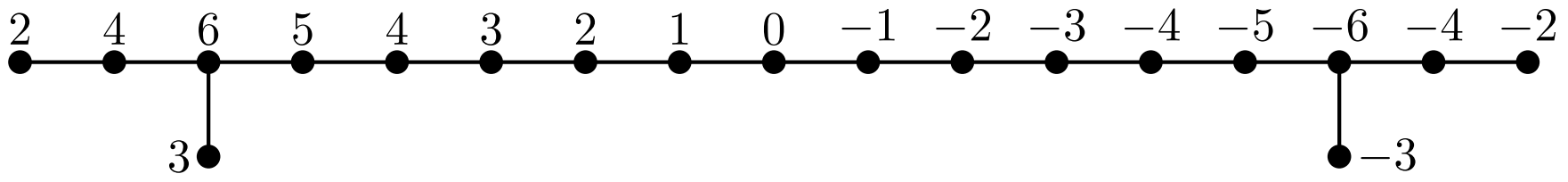
\begin{figure}
\begin{tikzpicture}[scale=0.7,thick]
\draw (2,0)--(2,1);
\draw (0,1)--(16,1);
\draw(14,0)--(14,1);

\node[above] at (0,1) {$2$};
\node[above] at (1,1) {$4$};
\node[above] at (2,1) {$6$};
\node[above] at (3,1) {$5$};
\node[above] at (4,1) {$4$};
\node[above] at (5,1) {$3$};
\node[above] at (6,1) {$2$};
\node[above] at (7,1) {$1$};
\node[above] at (8,1) {$0$};
\node[above] at (9,1) {$-1$};
\node[above] at (10,1) {$-2$};
\node[above] at (11,1) {$-3$};
\node[above] at (12,1) {$-4$};
\node[above] at (13,1) {$-5$};
\node[above] at (14,1) {$-6$};
\node[above] at (15,1) {$-4$};
\node[above] at (16,1) {$-2$};

\node[left] at (2,0) {$3$};
\node[right] at (14,0) {$-3$};

\filldraw[black]
(2,0) circle (3pt)
(14,0) circle (3pt)
(1,1) circle (3pt)
(0,1) circle (3pt)
(2,1) circle (3pt)
(3,1) circle (3pt)
(4,1) circle (3pt)
(5,1) circle (3pt)
(6,1) circle (3pt)
(7,1) circle (3pt)
(8,1) circle (3pt)
(9,1) circle (3pt)
(10,1) circle (3pt)
(11,1) circle (3pt)
(12,1) circle (3pt)
(13,1) circle (3pt)
(14,1) circle (3pt)
(15,1) circle (3pt)
(16,1) circle (3pt)
;

\end{tikzpicture}
\caption{The class $\xi$.}
\label{fig:Kgraph}
\end{figure}

Now fix labelling so that $K_{V_1}^2 \leq K_{V_2}^2$. Compute classes in $H^2(V_i;\bZ)$ as above; denote the classes in  $H^2(V_1;\bZ)$ by $f_j$ and those in $H^2(V_2;\bZ)$ by $f_j'$. Then define classes $r_i \in \xi^{\perp} \subset H^2(V_1;\bZ) \oplus H^2(V_2;\bZ)$, for $0 \leq i \leq 18$, as follows:
\begin{itemize}
\item If $K_{V_2}^2 = 9$, set $r_i = f_i$ for $0 \leq i \leq 17$ and set $r_{18} = f_{16} +2f_{17} +3f_{18} + f_0'$.
\item If $V_2 \cong \mathrm{Bl}_1(\bP^2) = \bF_1$, set $r_i = f_i$ for $0 \leq i \leq 16$, set $r_{17} =  f_{17} + f_1'$, and set $r_{18} = f_{16} +2f_{17}+ f_0'$.
\item If $V_2 \cong \bP^1 \times \bP^1$ or $V_2 \cong \bF_2$, set $r_i = f_i$ for $0 \leq i \leq 15$, set $r_{16} =  f_{16} + 2f_{17} + f_2'$, set $r_{17} = f_1'$, and set $r_{18} = f_{16}$.
\item  If $K_{V_2}^2 = 7$, set  $r_i = f_i$ for $0 \leq i \leq 15$, set $r_{16} = f_{16} + f_2'$, set $r_{17} = f_1'$, and set $r_{18} =f_{16} + f'_0$.
\item If $K_{V_2}^2 = d \leq 6$, set $r_i = f_i$ if $0 \leq i < 9+d$, set $r_{9+d} = f_{9+d}+f'_{9-d}$, and set $r_i = f'_{18-i}$ for $9+d < i \leq 18$.
\end{itemize}
Note that the second and third items above deal with the case that $K_{V_2}^2=8$. By construction, the classes $r_i$ all have self-intersection $(-2)$ and have the dual graph shown in Figure \ref{fig:dualgraph}. A straightforward linear algebra computation shows that the class $\xi$ is given by the linear combination of these classes shown in Figure \ref{fig:Kgraph}. 

The classes $\varphi_X(r_i)$, for $0 \leq i \leq 18$, form a set of simple roots in $I^{\perp}/I \cong H \oplus E_8 \oplus E_8$, which generate it as a $\bZ$-module. It follows that the $r_i$ form a $\bZ$-basis for $\xi^{\perp}$.

With this description, the lattices $\varphi_X(K_{V_i}^{\perp})$ are as follows. Here we take $E_3 = A_2 \oplus A_1$, $E_4 = A_4$, $E_5 = D_5$, and extend the $E$-series of lattices to $n \geq 9$ in the obvious way.
\begin{itemize}
\item If $K_{V_2}^2 = 9$, then $\varphi_X(K_{V_1}^{\perp}) = \langle r_0,\ldots,r_{17}\rangle \cong E_{18}$ and  $\varphi_X(K_{V_2}^{\perp}) = \{0\}$.
\item If $V_2 \cong \mathrm{Bl}_1(\bP^2) = \bF_1$, then $\varphi_X(K_{V_1}^{\perp}) = \langle r_0,\ldots,r_{16}\rangle \cong E_{17}$ and  $\varphi_X(K_{V_2}^{\perp}) = \langle r_{16} + 2r_{17} + r_{18} \rangle \cong \langle -8 \rangle$.
\item If $V_2 \cong \bP^1 \times \bP^1$ or $V_2 \cong \bF_2$, then $\varphi_X(K_{V_1}^{\perp}) = \langle r_0,\ldots,r_{15},r_{18}\rangle \cong E_{17}$ and  $\varphi_X(K_{V_2}^{\perp}) = \langle r_{17} \rangle \cong A_1$.
\item  If $K_{V_2}^2 = 7$, then $\varphi_X(K_{V_1}^{\perp}) = \langle r_0,\ldots,r_{15}\rangle \cong E_{16}$ and  $\varphi_X(K_{V_2}^{\perp}) = \langle r_{16}-r_{18}, r_{17} \rangle \cong \begin{psmallmatrix} -4 & 1 \\ 1 & -2 \end{psmallmatrix}$.
\item If $K_{V_2}^2 = d \leq 6$, then $\varphi_X(K_{V_1}^{\perp}) = \langle r_0,\ldots,r_{8+d}\rangle \cong E_{9+d}$ and  $\varphi_X(K_{V_2}^{\perp}) = \langle r_{10+d},\ldots,r_{18} \rangle \cong E_{9-d}$. 
\end{itemize}
Note once again that the second and third items above deal with the case that $K_{V_2}^2=8$. 
Therefore, given an $L$-polarized Tyurin degeneration of K3 surfaces, to describe the three sublattices $L$ and $\varphi_X(K_{V_i}^{\perp})$ of $I^{\perp}/I$ (up to isometries of $I^{\perp}/I$), it suffices to give the isometry types of the lattices $\varphi_X(K_{V_i}^{\perp})$ as abstract lattices, along with a set of generators for $L$ in terms of the simple roots $r_i \in I^{\perp}/I$ defined above. To denote such generators, we use the following shorthand notation: a class
\[v:= \sum_{i=0}^{18}a_ir_i \in I^{\perp}/I\]
will be denoted by the vector notation
\begin{equation} \label{eq:vectornotation} v = (a_0; a_1,a_2,a_3,a_4,a_5, a_6,a_7,a_8,a_9,a_{10},a_{11},a_{12},a_{13},a_{14},a_{15},a_{16},a_{17};a_{18}).\end{equation}
Note that any such class is only defined modulo 
\[\xi := (3;2,4,6,5,4,3,2,1,0,-1,-2,-3,-4,-5,-6,-4,-2;-3).\]
Having established these conventions, we have the following result.

\begin{proposition}\label{prop:degenerations}
Let $L$ denote one of the lattices $\langle 2 \rangle$, $\langle 4 \rangle$, or $H$. Then all very general $L$-polarized Tyurin degenerations of K3 surfaces, along with the corresponding sublattices $L$ and $\varphi_X(K_{V_i}^{\perp})$ of $I^{\perp}/I$, are listed in Table \ref{tab:degenerations}. All of these degenerations are doubly admissible.
\end{proposition}

\begin{proof}
The classification of $\langle 2 \rangle$-polarized Tyurin degenerations is essentially due to Friedman \cite[Section 5]{npgttk3s}; see also \cite{dk3sd2} and \cite[Section 2]{jonesthesis} for a more modern treatment. Descriptions of the restrictions of the polarizing class to the degenerate fibre are computed in \cite[Section 2]{jonesthesis}; the latter is given explicitly in terms of the classes $l$ and $e_i$ on each component, which can be rewritten in terms of the $r_i$ by simple linear algebra.

The computation for $\langle 4 \rangle$-polarized Tyurin degenerations is analogous, using the classification of $\langle 4 \rangle$-polarized Tyurin degenerations from \cite{t2dk3sd4}.

In the $H$-polarized case, we follow the approach from \cite{t2dk3sd4,jonesthesis}. Note that a general $H$-polarized K3 surface may be realised as a double cover of the Hirzebruch surface $\bF_4$ branched over the $(-4)$-section $s$ and a divisor $B$ in the linear system $|3s + 12f|$, where $f$ is the class of a fibre. There are two Type II boundary components in the moduli space of $H$-polarized K3 surfaces \cite[Section 7F]{cmek3sspt}, corresponding to two rank $2$ isotropic sublattices $I \subset H^{\perp}$, which we label by $D_{16}$ and $E_8E_8$. There are three easily-observed Tyurin degenerations:

\begin{itemize}[leftmargin=*]
\item \textbf{$D_{16}$ case.} Degenerate $B$ to the union of a double curve in the linear system $|s+4f|$ and a single curve in the linear system $|s+4f|$. The double cover is non-normal over the double curve; blowing this up gives rise to a Tyurin degeneration. The strict transform of the original central fibre is the Hirzebruch surface $V_2 := \bF_2$ and the exceptional component is $V_1 := \mathrm{Bl}_{17}\bP^2$.

\item \textbf{Asymmetric $E_8E_8$ case.} Degenerate $B$ to a curve with a consecutive triple point. The double cover has an elliptic singularity of type $\widetilde{E}_8$ over this point. Performing a weighted blow-up of this point (with weights $(1,1,2,3)$) gives rise to a Tyurin degeneration. The strict transform of the original central fibre is $V_1:= \mathrm{Bl}_{10}\bP^2$  and the exceptional component is $V_2:= \mathrm{dP}_1 = \mathrm{Bl}_{8}\bP^2$.

\item \textbf{Symmetric $E_8E_8$ case.} Degenerate $\bF_4$ to a union $\bF_2 \cup \bF_2$, glued along a fibre of the rulings. The $(-4)$-section degenerates to the union of the two $(-2)$-sections and the branch curve $B$ degenerates to a curve in the linear system $|3s+6f|$ in each component. The double cover is a Tyurin degeneration with both components isomorphic to  $\mathrm{Bl}_{9}\bP^2$.
\end{itemize}

Under the double cover, the classes $s$ and $f$ lift to the classes of a section $S$ and a fibre $F$ of the elliptic fibration on a general $H$-polarized K3 surface; it is straightforward to check that the flat limits $S_0$ and $F_0$ of these two divisors on the central fibre give the stated classes in Table \ref{tab:degenerations}. Note that in the symmetric $E_8E_8$ case we have a choice of which surface to label as $V_1$ and $V_2$; we choose the labelling so that $F_0$ lies in $V_2$.

It remains to prove that there are no further Tyurin degenerations. Following the approach of \cite[Section 5]{t2dk3sd4}, it suffices to compute the wall-and-chamber decomposition on the set of effective lifted polarizations.

Let $\calX \to \Delta$ denote one of the $H$-polarized Tyurin degenerations of K3 surfaces described above. Assume that $\NS(X) \cong H$ for a general fibre $X$ in $\calX$. Then the ample cone of $X$ is bounded by the rays spanned by the classes of $F$ and $S+2F$.

Let $A \subset \calX$ be an effective, flat, and relatively nef and big Cartier divisor, as in the definition of an $H$-polarized Tyurin degeneration (Definition \ref{def:tyurin}), and let $A_0$ be its restriction to $X_0$. The \emph{effective lifts} of $A$ are the set of divisors $A_0 + m \xi$ which are effective on $X_0$, where $m \in \bZ$. The \emph{set of effective lifted polarizations} is the set of all effective lifts.

The restriction to $X_0$ of any effective lift is linearly equivalent to a divisor of the form $aS_0 + bF_0 + c \xi$, for some $a,b,c \in \bZ$ and note that each of $S_0$, $F_0$, and $\xi$ is Cartier. Using the description of the divisors $S_0$ and $F_0$  from Table \ref{tab:degenerations} and the description of the ample cone above, we may compute the set of effective lifted polarizations in each case. We note first that the ample cone in $H$ gives the inequalities $a \geq 0$ and $b -2a \geq 0$ in both cases. 
\begin{itemize}[leftmargin=*]
\item \textbf{$D_{16}$ case.} The intersection of $F_0$ with each component of $X_0$ is an irreducible effective divisor of square $0$, so must intersect any other effective divisor non-negatively; this gives the further inequalities $c \geq 0$ and $a - 2c \geq 0$. These four inequalities describe the necessary conditions for an integral class of the above form to be effective. To show these conditions are sufficient, we consider the slice $a+b+c=k$, where $k\in \mathbb{R}$ and demonstrate that classes on the boundary rays are effective. By the convexity of the effective cone, this will provide the necessary and sufficient conditions for a lifted polarization to be effective.

For $k>0$, one can see that these conditions give rise to a triangle in the slice $a+b+c = k$ whose vertices are $(a,b) = (0,k),~ \left(\frac{2k}{7},\frac{4k}{7}\right),~ \left(\frac{k}{3},\frac{2k}{3}\right)$; we check that the $\mathbb{R}$-divisors corresponding to these vertices are effective. This is straightforward: the first and third vertices give effective sums of $S_0$ and $F_0$, and the second vertex gives the effective divisor $\frac{4k}{7}F_0|_{V_1} + \frac{k}{7}C|_{V_1}$ on $V_1$ and the zero divisor on $V_2$. Note that when $k=0$, the only class we obtain corresponds to the zero divisor and when $k<0$, the conditions describe the empty set.

\item \textbf{$E_8E_8$ case.} We present the analysis in the symmetric $E_8E_8$ case; the asymmetric case gives rise to the same picture. Each component $V_i$ of $X_0$ is a rational elliptic surface and the sum of a $(-1)$-section and a fibre of the elliptic fibration on $V_i$ is an effective divisor of square $1$, which intersects any other effective divisor non-negatively. This gives the further inequalities $c \geq 0$ and $b \geq c$. We proceed as in the $D_{16}$ case and consider the slice $a+b+c=k$.

For $k>0$, one can see that the conditions give rise to a quadrilateral in the slice $a+b+c = k$ whose vertices are $(a,b) = (0,k),~\left(0,\frac{k}{2}\right), ~\left(\frac{k}{5},\frac{2k}{5}\right),~\left(\frac{k}{3},\frac{2k}{3}\right)$. The first two vertices give rise to multiples of the effective anticanonical divisors $C|_{V_2}$ and $C|_{V_1}$, respectively. The fourth vertex gives $\frac{k}{3}S_0 + \frac{2k}{3}F_0 \sim \frac{k}{3}S_0|_{V_1} + \frac{2k}{3}C|_{V_2}$, which is effective, and the third gives the analogous effective $\bR$-divisor $\frac{k}{5}S_0|_{V_2} + \frac{2k}{5}C|_{V_1}$, with the labels of the components swapped. As before, this shows that the above conditions are necessary and sufficient conditions for effective lifted polarizations. For $k=0$, we obtain only the zero divisor and for $k<0$ the inequalities describe the empty set.
\end{itemize}

\noindent Within the set of effective lifted polarizations is a subset where the effective lifts are nef on $X_0$. For effective lifted polarizations outside this set, one may flip negative curves between components following the algorithm in the proof of \cite[Theorem 1]{bgd4} to obtain another Tyurin degeneration for which the effective lifted polarization becomes nef. This process endows the set of effective lifted polarizations with a wall-and-chamber structure, in which each chamber corresponds to a Tyurin degeneration and crossing walls corresponds to flipping curves. The argument from the proof of \cite[Theorem 2]{t2dk3sd4} shows that all Tyurin degenerations corresponding to a given Type II boundary component may be obtained in this way. In particular, we exploit the following facts: the nef cone is convex; intersection multiplicities are constant in linear systems; and the only curves which an effective reducible divisor may intersect negatively are its own components, which must necessarily have negative self-intersection.

Performing this analysis for the $H$-polarized degenerations above, we find that the divisors corresponding to the vertices given in the $D_{16}$ case are nef, and hence there is only one chamber in the set of effective lifted polarizations in the $D_{16}$ case. Therefore there is only one Tyurin degeneration in this case.

In the (symmetric) $E_8E_8$ case, the effective lifted polarization intersects a $(-1)$-section of the elliptic fibration on $V_1$ negatively if  $c < a $  and intersects a $(-1)$-section of the elliptic fibration on $V_2$ negatively if $b < a + c$. Flipping either of these sections gives the asymmetric $E_8E_8$ case (possibly with the labels of $V_1$ and $V_2$ switched). Working in a slice $a+b+c=k$, the two inequalities above can be reformulated to say that the conditions $b\geq k/2$ and $2a+b\geq k$ are necessary in order for a divisor to be nef in the symmetric model. This gives rise to a triangle with vertices $(0,k),~\left(0,\frac{k}{2}\right),~ \left(\frac{k}{4},\frac{k}{2}\right)$ and it is easy to check using the above well-known facts that all three of these vertices are nef, therefore this determines the nef cone for the symmetric model. Note that the first vertex intersects the $(-1)$-section on $V_2$ trivially, the second vertex intersects the $(-1)$-section on $V_1$ trivially, and the third vertex intersects both of these $(-1)$-sections trivially. 

Finally, it is easy to check that the only curve which the divisor associated to the vertex $\left(\frac{k}{5},\frac{2k}{5}\right)$ intersects negatively is the $(-1)$-section on $V_2$, hence after flipping this curve to $V_1$, the divisor corresponding to this vertex is a nef divisor on one version of the asymmetric model. Then all the vertices $\left(0,\frac{k}{2}\right),~\left(\frac{k}{5},\frac{2k}{5}\right),~ \left(\frac{k}{4},\frac{k}{2}\right)$ correspond to nef divisors on the asymmetric model and therefore the triangle defined by these points is the nef cone for this version of the asymmetric model. Analogous calculations performed on the vertex $\left(\frac{k}{3},\frac{2k}{3}\right)$ complete the analysis. 

The sets of effective lifted polarizations are depicted in Figure \ref{fig:cones}. In these pictures we depict the slices $a+b+c = 5$, with axes given by $a$ and $b$. The nodes demarcate points where $(a,b)$ are integers and the shaded region is the set of effective lifted polarizations. In the $D_{16}$ case there is only a single chamber. In the $E_8E_8$ case, the dark grey chamber corresponds to divisors which are nef in the symmetric case and the two light grey chambers correspond to divisors which are nef in the asymmetric cases.

\begin{figure}
\begin{tikzpicture}
\filldraw[lightgray] (0,5)--(1.4286,2.8571)--(1.6667,3.3333)--(0,5);

\draw[->] (0,0)--(3.5,0);
\draw[->] (0,0)--(0,5.5);
\node[right] at (3.5,0) {$a$};
\node[above] at (0,5.5) {$b$};
\draw (0,0)--(2.75,5.5);
\draw (0,5)--(3.5,1.5);
\draw (0,5)--(3.3333,0);

\filldraw (0,0) circle (2pt);
\filldraw (1,0) circle (2pt);
\filldraw (2,0) circle (2pt);
\filldraw (3,0) circle (2pt);
\filldraw (0,1) circle (2pt);
\filldraw (1,1) circle (2pt);
\filldraw (2,1) circle (2pt);
\filldraw (3,1) circle (2pt);
\filldraw (0,2) circle (2pt);
\filldraw (1,2) circle (2pt);
\filldraw (2,2) circle (2pt);
\filldraw (3,2) circle (2pt);
\filldraw (0,3) circle (2pt);
\filldraw (1,3) circle (2pt);
\filldraw (2,3) circle (2pt);
\filldraw (3,3) circle (2pt);
\filldraw (0,4) circle (2pt);
\filldraw (1,4) circle (2pt);
\filldraw (2,4) circle (2pt);
\filldraw (3,4) circle (2pt);
\filldraw (0,5) circle (2pt);
\filldraw (1,5) circle (2pt);
\filldraw (2,5) circle (2pt);
\filldraw (3,5) circle (2pt);

\node at (1.75,-0.5) {$D_{16}$ case};

\filldraw[lightgray] (6,5)--(6,2.5)--(7,2)--(7.6667,3.3333)--(6,5);
\filldraw[gray] (6,5)--(6,2.5)--(7.25,2.5)--(6,5);

\draw[->] (6,0)--(9.5,0);
\draw[->] (6,0)--(6,5.5);
\node[right] at (9.5,0) {$a$};
\node[above] at (6,5.5) {$b$};
\draw (6,0)--(8.75,5.5);
\draw (6,5)--(9.5,1.5);
\draw (6,2.5)--(9.5,0.75);
\draw (6,2.5)--(9.5,2.5);
\draw (6,5)--(8.5,0);

\filldraw (6,0) circle (2pt);
\filldraw (7,0) circle (2pt);
\filldraw (8,0) circle (2pt);
\filldraw (9,0) circle (2pt);
\filldraw (6,1) circle (2pt);
\filldraw (7,1) circle (2pt);
\filldraw (8,1) circle (2pt);
\filldraw (9,1) circle (2pt);
\filldraw (6,2) circle (2pt);
\filldraw (7,2) circle (2pt);
\filldraw (8,2) circle (2pt);
\filldraw (9,2) circle (2pt);
\filldraw (6,3) circle (2pt);
\filldraw (7,3) circle (2pt);
\filldraw (8,3) circle (2pt);
\filldraw (9,3) circle (2pt);
\filldraw (6,4) circle (2pt);
\filldraw (7,4) circle (2pt);
\filldraw (8,4) circle (2pt);
\filldraw (9,4) circle (2pt);
\filldraw (6,5) circle (2pt);
\filldraw (7,5) circle (2pt);
\filldraw (8,5) circle (2pt);
\filldraw (9,5) circle (2pt);

\node at (7.75,-0.5) {$E_8E_8$ case};
\end{tikzpicture}
\caption{Wall and chamber decompositions of the sets of effective lifted polarizations} \label{fig:cones}
\end{figure}
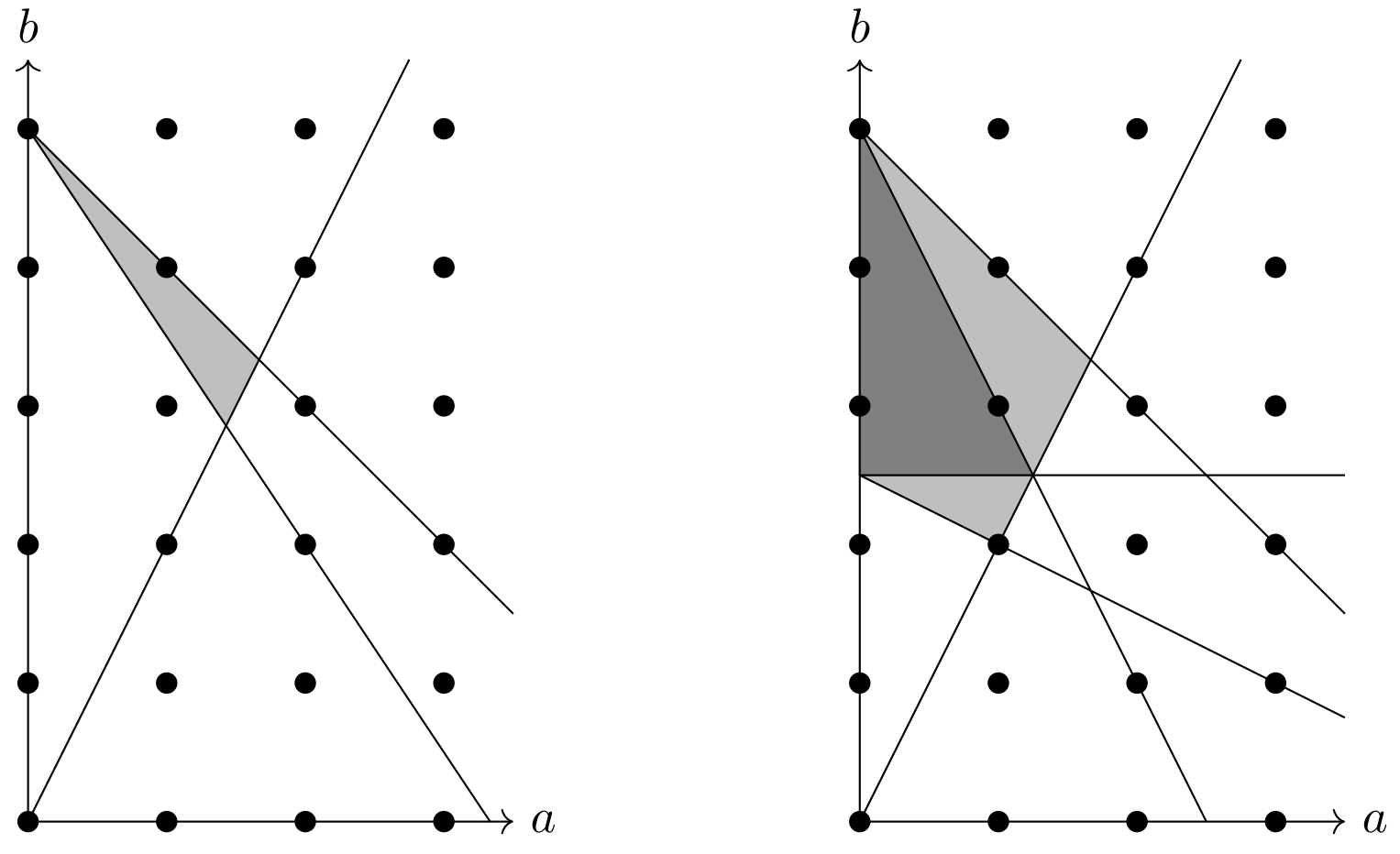

Finally, we check that these degenerations are doubly admissible. Indeed, the two rank two isotropic sublattices $I \subset H^{\perp} \cong H \oplus H \oplus E_8 \oplus E_8$ correspond to the two embeddings of $H \oplus H$ into $H \oplus H \oplus E_8 \oplus E_8$: the obvious one with orthogonal complement $E_8 \oplus E_8$ and a second with orthogonal complement given by an index $2$ overlattice of $D_{16}$. The doubly admissible condition then follows from Remark \ref{rem:doublyadmissible}.
\end{proof}

\subsection{Mirror symmetry}\label{sec:mirrorsymmetry}

In \cite[Conjecture 5.10]{dfk3slpms}, it is conjectured that mirror symmetry exchanges Tyurin degenerations of a K3 surface and allowable splittings of its mirror K3 surface.  A special case of this conjecture can be stated as follows (c.f. \cite[Proposition 5.3 and Remark 5.12]{dfk3slpms}).

\begin{conjecture}\label{conj:dht}
Let $\calX \to \Delta$ be a very general doubly admissible $L$-polarized Tyurin degeneration of K3 surfaces. Let $\check{L}$ be the mirror lattice to $L$, in the sense of \cite[Section 6]{mslpk3s}. Then a very general $\check{L}$-quasipolarized K3 surface $Y$ admits an elliptic fibration $\pi\colon Y \to \bP^1$ along with an allowable loop $\gamma \subset \bP^1$, such that the $\check{L}$-quasipolarization is doubly admissible and compatible with $\pi$ and $\gamma$, and there is an isometry $\psi\colon J^{\perp}/J \to I^{\perp}/I$ such that $\psi(\Gamma) = L^{\perp}$ and $\psi(\varphi_Y([r_i(b)]^{\perp})) = \varphi_X(K_{V_i}^\perp)$ for each $i \in \{1,2\}$.
\end{conjecture}

We prove three special cases of this conjecture.

\begin{theorem}\label{thm:mirror}
Let $\calX \to \Delta$ be a very general $\langle 2 \rangle$-polarized (resp. $\langle 4 \rangle$-polarized, $H$-polarized) Tyurin degeneration of K3 surfaces. Then there exists a very general $M_1$-quasi\-polarized (resp. $M_2$-quasi\-polarized, $M$-quasi\-polarized)  K3 surface $Y$ along with an elliptic fibration $\pi\colon Y \to \bP^1$ and an allowable loop $\gamma \subset \bP^1$, such that the $\check{L}$-quasipolarization is doubly admissible and compatible with $\pi$ and $\gamma$, and there is an isometry $\psi\colon J^{\perp}/J \to I^{\perp}/I$ such that $\psi(\Gamma) = L^{\perp}$ and $\psi(\varphi_Y([r_i(b)]^{\perp})) = \varphi_X(K_{V_i}^\perp)$ for each $i \in \{1,2\}$.
\end{theorem}

\begin{proof}
The matching between Tyurin degenerations and fibrations is simply given by comparing rows in Tables \ref{tab:splittingdata} and \ref{tab:degenerations} (noting that mirror symmetry predicts that $e(Y_i) - 12 = -K_{V_i}^2$). The sublattices $\Gamma$ and $\varphi_Y([r_i(b)]^{\perp})$ of $J^{\perp}/J \cong H \oplus E_8 \oplus E_8$ are computed explicitly in Section \ref{sec:explicitsplitting} and the sublattices $L$ and $\varphi_X(K_{V_i}^\perp)$ of $I^{\perp}/I \cong H \oplus E_8 \oplus E_8$ are computed explicitly in Proposition \ref{prop:degenerations}. To prove the final part it suffices to find an explicit isometry of $H \oplus E_8 \oplus E_8$ matching $\Gamma$ with $L^{\perp}$ and matching $\varphi_Y([r_i(b)]^{\perp})$ with $\varphi_X(K_{V_i}^\perp)$, for $i\in\{1,2\}$, in each case. 

Our approach to this is quite straightforward. In the cases under consideration, the sublattices $\Gamma \subset J^{\perp}/J$ and $L^{\perp} \subset I^{\perp}/I$ are large (i.e. corank $1$ or $2$) and negative definite. The intersections $\Gamma \cap \varphi_Y([r_i(b)]^{\perp})$ and  $L^{\perp} \cap \varphi_X(K_{V_i}^\perp)$, for each $i \in \{1,2\}$, are also negative definite. There are thus a finite number of isometries between  $\Gamma \cap \varphi_Y([r_i(b)]^{\perp})$ and  $L^{\perp} \cap \varphi_X(K_{V_i}^\perp)$;  we performed a brute-force search over these isometries, to find ones which admit lifts to isometries between $\Gamma$ and $L^{\perp}$, and thence between $J^{\perp}/J$ and $I^{\perp}/I$, which have the required properties. SageMath  worksheets replicating these computations are available at \cite{DHTcode}.
\end{proof}

As an application of this result, we perform an investigation into the properties of the \emph{coupling group}. The coupling group is a somewhat mysterious object which relates the lattice polarization on the Tyurin degeneration (resp. elliptically fibred K3 surface) to the induced lattice polarizations on the components of the degenerate fibre (resp. splitting). It is defined as follows; compare \cite[Definition 2.35]{dfk3slpms}.

\begin{definition}
Let $\calX \to \Delta$ be an $L$-polarized Tyurin degeneration of K3 surfaces. Define $L_i := L \cap \varphi_X(K_{V_i}^{\perp})$. Then the \emph{coupling group} is defined to be
\[Q(L) := L / \langle L_1, L_2 \rangle.\]
Similarly, let  $Y$ be an $\check{L}$-quasipolarized K3 surface admitting an elliptic fibration $\pi\colon Y \to \bP^1$, along with an allowable loop $\gamma \subset \bP^1$, such that the $\check{L}$-quasipolarization is compatible with $\pi$ and $\gamma$. Define $\Gamma_i := \Gamma \cap \varphi_Y([r_i(b)]^{\perp})$. Then the \emph{coupling group} is defined to be
\[Q(\Gamma) := \Gamma / \langle \Gamma_1, \Gamma_2 \rangle.\]
The coupling groups $Q(L)$ and $Q(\Gamma)$ are said to be \emph{torsion} if they are of rank 0.
\end{definition}

Various properties of the coupling group are proved in \cite{dfk3slpms}; the following corollaries elaborate on these properties. Our first result shows that the coupling group is not an invariant of the lattice-theoretic data defining the degeneration (resp. fibration): specifically, it is not uniquely determined by the polarizing lattice $L$ and the isotropic sublattice $I \subset L^{\perp}$ (resp. $\check{L}$ and the choice of fibration $\pi\colon Y \to \bP^1$), but depends upon the choice of birational model (resp. splitting loop) in a nontrivial way.

\begin{corollary}\label{cor:coupling1} The coupling group $Q(L)$ (resp. $Q(\Gamma)$) depends upon the choice of birational model for the Tyurin degeneration (resp. splitting loop).
\end{corollary}

\begin{proof} We may verify these statements by direct computation of coupling groups in some of the examples from Propositions \ref{prop:splittings} and \ref{prop:degenerations}. The two birational models for a $\langle 4 \rangle$-polarized Tyurin degeneration of type $E_8E_8\langle -4 \rangle$ have coupling groups $Q(L) = 0$ in the asymmetric case with $(K_{V_1}^2,K_{V_2}^2) = (-1,1)$ and $Q(L) = \bZ$ in the symmetric case with $(K_{V_1}^2,K_{V_2}^2) = (0,0)$. A similar statement holds for the coupling groups $Q(\Gamma)$ of the mirror $M_2$-quasipolarized elliptic fibrations from Theorem \ref{thm:mirror}.
\end{proof}

Our second result is a strengthening of \cite[Proposition 4.27]{dfk3slpms}, which describes the relationship between $Q(\Gamma)$ and the Mordell-Weil group of $\pi\colon Y \to \bP^1$.

\begin{corollary}\label{cor:coupling2}
 There is a surjective map $\MW(Y,\pi) \to Q(\Gamma)$, but these two groups do not need to have the same rank. In particular, $\MW(Y,\pi)$ may have strictly positive rank, even if $Q(\Gamma)$ is torsion. 
\end{corollary}
\begin{proof}
 The map $\MW(Y,\pi) \to Q(\Gamma)$ is surjective by \cite[Proposition 4.27]{dfk3slpms}, but note that the asymmetric splitting $\rmI\rmI^*\, 3\rmI_1 \, |\, \rmI\rmI^*\, \rmI_1$ of the $E_8E_8\langle -4 \rangle$ fibration on an $M_2$-quasipolarized K3 surface has $Q(\Gamma) = 0$ and $\MW(Y,\pi) = \bZ$, so these groups need not have the same rank. 
\end{proof}

Our final result is a strengthening of \cite[Lemma 5.7]{dfk3slpms}, which explores the relationship between $Q(L)$ and $Q(\Gamma)$ for a mirror pair consisting of an $L$-polarized Tyurin degeneration and an $\check{L}$-quasipolarized elliptic fibration.

\begin{corollary}\label{cor:coupling3}
Let $\calX \to \Delta$ be a very general $L$-polarized Tyurin degeneration and let $\pi\colon Y \to \bP^1$ be an $\check{L}$-quasipolarized elliptically fibred K3 surface with a splitting loop $\gamma \subset \bP^1$ which is mirror to it in the sense of Conjecture \ref{conj:dht}.
\begin{enumerate}
\item If $e(Y_1) - 12 = -K_{V_1}^2 \neq 0$, then $Q(L)$ is torsion if and only if $Q(\Gamma)$ is torsion.
\item If $e(Y_1) - 12 = -K_{V_1}^2 = 0$, then one of $Q(L)$ and $Q(\Gamma)$ may be torsion and the other may be non-torsion.
\item If $e(Y_1) - 12 = -K_{V_1}^2 \neq 0$ and $Q(L)$ and $Q(\Gamma)$ are both torsion, we may nevertheless have $Q(L) \neq Q(\Gamma)$. In this case one of $Q(L)$ and $Q(\Gamma)$ may be trivial whilst the other is nontrivial.
\item If $e(Y_1) - 12 = -K_{V_1}^2 \neq 0$ and $Q(L)$ and $Q(\Gamma)$ are both non-torsion, we may nevertheless have $Q(L) \neq Q(\Gamma)$.
\end{enumerate}
\end{corollary}
\begin{proof}
We verify these statements by direct computation of coupling groups in some of the examples from Theorem \ref{thm:mirror}.
\begin{enumerate}
\item This is precisely \cite[Lemma 5.7]{dfk3slpms}.
\item The symmetric $H$-polarized Tyurin degeneration of type $E_8E_8$ has $K_{V_1}^2= 0$. Its coupling group is $Q(L) = \bZ$ and its mirror fibration has coupling group $Q(\Gamma) = 0$. 
\item The $H$-polarized Tyurin degeneration of type $D_{16}$ has $K_{V_1}^2= -8$ and $Q(L) = \frac{\bZ}{4\bZ}$, whilst its mirror fibration has $Q(\Gamma) = \frac{\bZ}{2\bZ}$. The $\langle 4 \rangle$-polarized Tyurin degeneration of type $D_{17}$ has $K_{V_1}^2= -9$ and $Q(L) = \frac{\bZ}{3\bZ}$, whilst its mirror fibration has $Q(\Gamma) = 0$.
\item The $\langle 4 \rangle$-polarized Tyurin degeneration of type $D_{16}\langle -4 \rangle$ has $K_{V_1}^2= -8$ and $Q(L) = \bZ$, whilst its mirror fibration has $Q(\Gamma) = \bZ \oplus \frac{\bZ}{2\bZ}$.
\end{enumerate}
\end{proof}

\appendix
\section{Tables of elliptic fibrations} \label{app:elliptictables}

Tables \ref{tab:Mweierstrass}, \ref{tab:M1weierstrass}, and \ref{tab:M2weierstrass} contain the data of the elliptically fibred K3 surfaces studied in this paper. In these tables:
\begin{itemize}
\item ``Label'' is the label used to refer to each example in the text.
\item ``Fibres'' is the singular fibre type of a very general member of the family.
\item ``Weierstrass parameters'' are the functions $a_t(u)$ and $b_t(u)$ appearing in the Weierstrass form $y^2 = x^3 + a_t(u)x + b_t(u)$. Here $u$ is an affine parameter on the base $\bP^1$ and $t = (t_1,\ldots,t_n)$ are deformation parameters.
\item ``$\bQ$ point'' is a point $(t_1,\ldots,t_n)$ at which the functions $a_t(u)$ and $b_t(u)$ become polynomials in $\bQ[u]$. Substituting this point into the Weierstrass equation gives the surface upon which the algorithm from \cite{snahlpcesp1} was run.
\item ``MW'' is the Mordell-Weil lattice of a very general member of the family.
\end{itemize}

\begin{landscape}
\begingroup
\renewcommand{\arraystretch}{1.1}
\begin{table}
\begin{tabular}{|c|c|c|c|c|}
\hline
 \textbf{Label} & \textbf{Fibres} & \textbf{Weierstrass parameters}& \textbf{$\bQ$ point} & \textbf{MW}\\ 
\hline
 \multirow{2}*{$D_{16}$} &  \multirow{2}*{$\rmI_{12}^*\, 6\rmI_1$} & $a_t(u) = u^2(u^6 - \frac{1}{3}(1 - 3t_1u^2 - 2t_2u^3)^2)$ &  \multirow{2}*{$(-3,7)$} &  \multirow{2}*{$\frac{\bZ}{2 \bZ}$}  \\
  &  & $b_t(u) = \frac{1}{27}u^3(1 - 3t_1u^2 - 2t_2u^3)(2(1 - 3t_1u^2 - 2t_2u^3)^2 - 9u^6)$ & &  \\
\hline
 \multirow{2}*{$E_8E_8$} &  \multirow{2}*{$2\rmI\rmI^*\, 4\rmI_1$} & $a_t(u) = -3t_1u^4$ &  \multirow{2}*{$(-3,7)$} &  \multirow{2}*{$1$}  \\
  &  & $b_t(u) = -u^5(1-2t_2u+u^2)$ & &    \\
\hline
\end{tabular}
\caption{Elliptic fibrations on $M$-quasipolarized K3 surfaces.}
\label{tab:Mweierstrass}
\end{table}

\begin{table}
\begin{tabular}{|c|c|c|c|c|}
\hline
 \textbf{Label} & \textbf{Fibres} & \textbf{Weierstrass parameters}& \textbf{$\bQ$ point} & \textbf{MW}\\ 
\hline
 \multirow{2}*{$A_{17}$} &  \multirow{2}*{$\rmI_{18}\, 6\rmI_1$} & $a_t(u) = (t_1+t_2u+t_3u^2)\left(1-\frac{1}{48t_4^2}(t_1+t_2u+t_3u^2)^3\right)$ &  \multirow{2}*{$(3,7,-3,5)$} &  \multirow{2}*{$\frac{\bZ}{3 \bZ}$}  \\
 &  & $b_t(u) = \frac{1}{864t_4^3}(t_1+t_2u+t_3u^2)^6 - \frac{1}{12t_4}(t_1+t_2u+t_3u^2)^3+t_4$ & &   \\
\hline
 \multirow{3}*{$D_{16}A_1$} &  \multirow{3}*{$\rmI_{12}^*\,\rmI_2 \, 4\rmI_1$} & $a_t(u) = u^2(u^6 - \frac{1}{3}(1 - 3t_1u^2 - 2t_2u^3)^2)$ &  \multirow{3}*{$(4,7)$} &  \multirow{3}*{$\frac{\bZ}{2 \bZ}$}  \\
  &  & $b_t(u) = \frac{1}{27}u^3(1 - 3t_1u^2 - 2t_2u^3)(2(1 - 3t_1u^2 - 2t_2u^3)^2 - 9u^6)$ & &  \\
   & & $t_1^3 = (t_2 \pm 1)^2$ & & \\
\hline
 \multirow{2}*{$E_7D_{10}$} &  \multirow{2}*{$\rmI\rmI\rmI^*\, \rmI_{6}^*\,3\rmI_1$} & $a_t(u) = (t_3u^3 - 3(t_1+t_2u)^2)u^2$ &  \multirow{2}*{$(3,7,-3)$} &  \multirow{2}*{$\frac{\bZ}{2 \bZ}$} \\
  &  & $b_t(u) = u^3(t_1+t_2u)(2(t_1+t_2u)^2 - t_3u^3)$ & &   \\
\hline
 \multirow{3}*{$E_8E_8A_1$} &  \multirow{3}*{$2\rmI\rmI^*\, \rmI_2\, 2\rmI_1$} & $a_t(u) = -3t_1u^4$ &  \multirow{3}*{$(4,9)$} &  \multirow{3}*{$1$}  \\
  &  & $b_t(u) = -u^5(1-2t_2u+u^2)$ & &    \\
     & & $t_1^3 = (t_2 \pm 1)^2$ & & \\
\hline
\end{tabular}
\caption{Elliptic fibrations on $M_1$-quasipolarized K3 surfaces.}
\label{tab:M1weierstrass}
\end{table}

\begin{table}
\begin{tabular}{|c|c|c|c|c|}
\hline
 \textbf{Label} & \textbf{Fibres} & \textbf{Weierstrass parameters}& \textbf{$\bQ$ point} & \textbf{MW}\\ 
\hline
 \multirow{2}*{$A_{15}A_1A_1$} &  \multirow{2}*{$\rmI_{16}\, 2\rmI_2\, 4\rmI_1$} & $a_t(u) = -\frac{1}{3}(t_1^2u^4(1-u)^4   + 4t_1u^2(1-u)^2 + 1)$ &  \multirow{2}*{$3$} &  \multirow{2}*{$\frac{\bZ}{4 \bZ}$}  \\
 &  & $b_t(u) = \frac{1}{27}(t_1u^2(1-u)^2 + 2)(2t_1^2u^4(1-u)^4 + 8t_1u^2(1-u)^2 - 1)$ & &   \\
\hline
 \multirow{2}*{$D_{17}$} &  \multirow{2}*{$\rmI_{13}^*\,5\rmI_1$} & $a_t(u) = -\frac{27}{4}u^2(4 - 8t_1^3u^3 + 54t_1^4u^4 + 54t_1^5u^5 - 23t_1^6u^6)$ &  \multirow{2}*{$3$} &  \multirow{2}*{$1$} \\
  &  & $b_t(u) = \frac{27}{16}u^3(32 - 96t_1^3u^3 + 648t_1^4u^4 + 648t_1^5u^5- 228t_1^6u^6 - 648t_1^7u^7 + 1539t_1^8u^8 + 4666t_1^9u^9 )$ & &   \\
\hline
 \multirow{3}*{$D_{16}\langle -4 \rangle$} &  \multirow{3}*{$\rmI_{12}^*\, 6\rmI_1$} & $a_t(u) = u^2(u^6 - \frac{1}{3}(1 - 3t_1u^2 - 2t_2u^3)^2)$ &  \multirow{3}*{$(5,-1)$} &  \multirow{3}*{$\bZ \times \frac{\bZ}{2 \bZ}$}  \\
  &  & $b_t(u) = \frac{1}{27}u^3(1 - 3t_1u^2 - 2t_2u^3)(2(1 - 3t_1u^2 - 2t_2u^3)^2 - 9u^6)$ & &  \\
   & & $(t_1,t_2)$ satisfy Equation \eqref{eq:4polNL} & & \\
\hline
 \multirow{2}*{$D_{12}D_5$} &  \multirow{2}*{$\rmI_8^*\, \rmI_{1}^*\,3\rmI_1$} & $a_t(u) = u^2(1-\frac{1}{3}(2+t_1u+t_2u^2)^2)$ &  \multirow{2}*{$(3,7)$} &  \multirow{2}*{$\frac{\bZ}{2 \bZ}$} \\
  &  & $b_t(u) = \frac{1}{27}u^3(2+t_1u+t_2u^2)(2(2+t_1u+t_2u^2)^2 - 9)$ & &   \\
\hline
 \multirow{2}*{$D_8D_8 \langle -4 \rangle$} &  \multirow{2}*{$2\rmI_4^*\, 4\rmI_1$} & $a_t(u) =-\frac{1}{3}u^2(t_1^2 + 2t_1t_2u + t_2^2u^2 - t_1t_3u^2 + 2t_2t_3u^3 + t_3^2u^4)$ &  \multirow{2}*{$(5,9,-1)$} &  \multirow{2}*{$\bZ \times \frac{\bZ}{2 \bZ}$}  \\
 &  & $b_t(u) = \frac{1}{27}u^3(t_1 + t_2u + t_3u^2)(2t_1^2 + 4t_1t_2u + 2t_2^2u^2 - 5t_1t_3u^2 + 4t_2t_3u^3 + 2t_3^2u^4)$ & &   \\
\hline
 \multirow{3}*{$E_8E_8\langle -4 \rangle$} &  \multirow{3}*{$2\rmI\rmI^*\, 4\rmI_1$} & $a_t(u) = -3t_1u^4$ &  \multirow{3}*{$(5,-1)$} &  \multirow{3}*{$\bZ$}  \\
  &  & $b_t(u) = -u^5(1-2t_2u+u^2)$ & &    \\
     & & $(t_1,t_2)$ satisfy Equation \eqref{eq:4polNL} & & \\
\hline
 \multirow{2}*{$E_8D_9$} &  \multirow{2}*{$\rmI\rmI^*\, \rmI^*_5\, 3\rmI_1$} & $a_t(u) = -u^2 (3t_1^2 - 6t_1t_2u - t_3u^2)$ &  \multirow{2}*{$(3,7,-3)$} &  \multirow{2}*{$1$}  \\
 &  & $b_t(u) = \frac{1}{12t_1} u^3(24t_1^4 - 72t_1^3t_2u + 12t_1^2(3t_2^2 - t_3)u^2 + 12t_1t_2(t_2^2 + t_3)u^3 + (9t_2^4 + 6t_2^2t_3 + t_3^2)u^4)$ & &   \\
\hline
 \multirow{2}*{$E_7E_7A_3$} &  \multirow{2}*{$2\rmI\rmI\rmI^*\, \rmI_4\, 2\rmI_1$} & $a_t(u) = u^3(t_2  - 3t_1^2u - 2t_2u  + t_2u^2)$ &  \multirow{2}*{$(3,7)$} &  \multirow{2}*{$\frac{\bZ}{2 \bZ}$}  \\
 &  & $b_t(u) = t_1u^5(- t_2 + 2t_1^2u  + 2t_2u - t_2u^2)$ & &   \\
\hline
 \multirow{3}*{$E_6A_{11}$} &  \multirow{3}*{$\rmI\rmV^*\, \rmI_{12}\, 4\rmI_1$} & $a_t(u) = -\frac{1}{3}t_1(t_3u + 1)(t_1 + 3t_1t_3u + 3t_1t_3^2u^2 + t_1t_3^3u^3 - 6t_2u^4)$ &  \multirow{3}*{$(3,5,-3)$} &  \multirow{3}*{$\frac{\bZ}{3 \bZ}$} \\
  &  & $\begin{array}{l} b_t(u) = \frac{1}{27}t_1(2t_1^2 + 12t_1^2t_3u + 30t_1^2t_3^2u^2 + 40t_1^2t_3^3u^3 + (30t_1^2t_3^4 - 18t_1t_2)u^4 \\ \qquad \qquad+ (12t_1^2t_3^5 - 54t_1t_2t_3)u^5 + (2t_1^2t_3^6     - 54t_1t_2t_3^2)u^6    - 18t_1t_2t_3^3u^7 + 27t_2^2u^8)\end{array}$ & &   \\
\hline
\end{tabular}
\caption{Elliptic fibrations on $M_2$-quasipolarized K3 surfaces.}
\label{tab:M2weierstrass}
\end{table}
\endgroup
\end{landscape}

\section{Tables of mirror symmetry data} \label{app:mirrortables}

Table \ref{tab:splittingdata} contains the data for the computations of splittings of the elliptic fibrations from Section \ref{sec:splittings}. The fibrations are those given in Tables \ref{tab:Mweierstrass}, \ref{tab:M1weierstrass}, and \ref{tab:M2weierstrass}. In this table:
\begin{itemize}
\item ``Lattice'' is the polarizing lattice of the elliptically fibred K3 surface.
\item ``Label'' is the label used to refer to each example in the text.
\item ``Splittings'' describes the singular fibre types on either side of the allowable loop; the loop is represented by a vertical bar $|$. We label the pieces of the splitting so that $Y_1$ contains the fibres listed to the left of the bar and $Y_2$ contains those listed to the right.
\item ``Monodromy'' gives the values $(e(Y_1) - 12, e(Y_2)-12)$ in the monodromy matrix from Definition \ref{def:allowable}(1) on each side of the loop.
\item  ``$\varphi_Y([r_i(b)]^{\perp})$'' describes the pair of lattices $(\varphi_Y([r_1(b)]^{\perp}), \varphi_Y([r_2(b)]^{\perp}))$ as abstract lattices. Here we extend the $E_n$-series of lattices to $n \geq 9$ in the obvious way.
\end{itemize}

\begin{table}
\begin{tabular}{|c|c|c|c|c|}
\hline
\textbf{Lattice} & \textbf{Label} & \textbf{Splittings} & \textbf{Monodromy} & $\varphi_Y([r_i(b)]^{\perp})$\\ 
\hline
\hline
\multirow{3}*{$M$} & $D_{16}$ & $\rmI_{12}^*\,  2\rmI_1\, |\, 4\rmI_1$ & $(8, -8)$ & $(E_{17}, A_1)$\\
\cline{2-5}
&\multirow{2}*{$E_8E_8$} &  $\rmI\rmI^*\, 3\rmI_1 \, |\, \rmI\rmI^*\, \rmI_1$ & $(1,-1)$ & $(E_{10}, E_8)$ \\
& & $\rmI\rmI^*\, 2\rmI_1 \, |\, \rmI\rmI^*\, 2\rmI_1$ & $(0,0)$& $(E_{9}, E_9)$ \\
\hline
\hline
\multirow{4}*{$M_1$} & $A_{17}$ & $\rmI_{18}\, 3\rmI_1\, |\, 3\rmI_1$ & $(9, -9)$& $(E_{18}, 0)$ \\
\cline{2-5}
&$D_{16}A_1$ & $\rmI_{12}^*\, 2\rmI_1 \, | \, \rmI_2 \, 2\rmI_1$ & $(8,-8)$&  $(E_{17}, A_1)$\\
\cline{2-5}
& $E_7D_{10}$ & $\rmI\rmI\rmI^*\, \rmI_1\, |\, \rmI_{6}^*\,\rmI_1$ & $(2, -2)$& $(E_{11}, E_7)$ \\
\cline{2-5}
&$E_8E_8A_1$ & $\rmI\rmI^*\, \rmI_2\, \rmI_1 \, |\, \rmI\rmI^*\, \rmI_1$& $(1,-1)$& $(E_{10}, E_8)$ \\
\hline
\hline
\multirow{11}*{$M_2$} & $A_{15}A_1A_1$ & $\rmI_{16}\, \rmI_2\, 2\rmI_1 \,|\, \rmI_2\, 2\rmI_1$  & $(8, -8)$& $(E_{17}, A_1)$ \\
\cline{2-5}
&$D_{17}$ & $\rmI_{13}^*\,2\rmI_1\, | \, 3\rmI_1$ & $(9,-9)$& $(E_{18}, 0)$ \\
\cline{2-5}
&$D_{16}\langle -4 \rangle$ & $\rmI_{12}^*\,  2\rmI_1\, |\, 4\rmI_1$ & $(8, -8)$& $(E_{17}, A_1)$ \\
\cline{2-5}
&$D_{12}D_5$ & $\rmI_8^*\, 2\rmI_1\, |\, \rmI_{1}^*\,\rmI_1$& $(4,-4)$& $(E_{13}, D_5)$ \\
\cline{2-5}
&$D_8D_8 \langle -4 \rangle$ & $\rmI_4^*\, 2\rmI_1\, |\, \rmI_4^*\, 2\rmI_1$& $(0,0)$& $(E_{9}, E_9)$ \\
\cline{2-5}
&\multirow{2}*{$E_8E_8\langle -4 \rangle$} &  $\rmI\rmI^*\, 3\rmI_1 \, |\, \rmI\rmI^*\, \rmI_1$ & $(1,-1)$ & $(E_{10}, E_8)$ \\
& & $\rmI\rmI^*\, 2\rmI_1 \, |\, \rmI\rmI^*\, 2\rmI_1$ & $(0,0)$& $(E_{9}, E_9)$ \\
\cline{2-5}
&$E_8D_9$ & $\rmI^*_5\, 2\rmI_1\, | \, \rmI\rmI^*\, \rmI_1$ & $(1,-1)$& $(E_{10}, E_8)$ \\
\cline{2-5}
&$E_7E_7A_3$ & $\rmI\rmI\rmI^*\, \rmI_4\, \rmI_1\, | \, \rmI\rmI\rmI^*\, \rmI_1$ & $(2,-2)$ & $(E_{11}, E_7)$\\
\cline{2-5}
&\multirow{2}*{$E_6A_{11}$} &  $\rmI\rmV^*\, \rmI_{12}\, \rmI_1\, | \, 3\rmI_1$ & $(9,-9)$& $(E_{18}, 0)$  \\
& & $\rmI_{12}\, 3\rmI_1\, | \, \rmI\rmV^*\, \rmI_1$ & $(3,-3)$& $(E_{12}, E_{6})$ \\
\hline
\end{tabular}
\caption{Splitting data.}
\label{tab:splittingdata}
\end{table}

Table \ref{tab:degenerations} contains the data for the Tyurin degenerations considered in Section \ref{sec:tyurin}. In this table:
\begin{itemize}
\item ``Lattice'' is the polarizing lattice of the Tyurin degeneration.
\item ``Label'' is the label used to refer to each example in the text.
\item ``$V_1$'' and ``$V_2$'' are the two components of the degenerate fibre $X_0$.
\item ``$(K_{V_1}^2,K_{V_2}^2)$'' gives the self-intersection numbers of the canonical classes on each of $V_1$ and $V_2$.
\item  ``$\varphi_X(K_{V_i}^\perp)$'' describes the pair of lattices $(\varphi_X(K_{V_1}^\perp), \varphi_X(K_{V_2}^\perp))$ as abstract lattices. Here we extend the $E_n$-series of lattices to $n \geq 9$ in the obvious way.
\item  ``Generators of $L$'' describes generators of the lattice $L$ in terms of the simple roots $r_i \in I^{\perp}/I$. The notation is described in Equation \eqref{eq:vectornotation}.
\end{itemize}

\begin{landscape}
\begingroup
\renewcommand{\arraystretch}{1.1}
\begin{table}
\begin{tabular}{|c|c|c|c|c|c|c|}
\hline
\textbf{Lattice} & \textbf{Label} & $V_1$ & $V_2$ & $(K_{V_1}^2,K_{V_2}^2)$ & $\varphi_X(K_{V_i}^\perp)$ & \textbf{Generators of} $L$\\
\hline
\hline 
\multirow{6}*{$H$} & \multirow{2}*{$D_{16}$}   & \multirow{2}*{$\mathrm{Bl}_{17}\bP^2$}    & \multirow{2}*{$\bF_2$}  & \multirow{2}*{$(-8,8)$} & \multirow{2}*{$(E_{17},A_1)$}  & $(0;0,0,0,0,0,0,0,0,0,0,0,0,0,0,0,0,1;0)$ \\ 
&   &   &  &  &  & $(1;0,1,2,2,2,2,2,2,2,2,2,2,2,2,2,1,0;1)$ \\ 
\cline{2-7}
& \multirow{4}*{$E_8E_8$}  & \multirow{2}*{$\mathrm{Bl}_{10}\bP^2$}    & \multirow{2}*{$\mathrm{Bl}_{8}\bP^2$} & \multirow{2}*{$(-1,1)$} & \multirow{2}*{$(E_{10},E_8)$}  & $(0;0,0,0,0,0,0,0,0,1,0,0,0,0,0,0,0,0;0)$ \\ 
&   &   &  &  &  & $(3;2,4,6,5,4,3,2,1,0,0,0,0,0,0,0,0,0;0)$ \\ 
\cline{3-7}
& & \multirow{2}*{$\mathrm{Bl}_{9}\bP^2$}    & \multirow{2}*{$\mathrm{Bl}_{9}\bP^2$}  & \multirow{2}*{$(0,0)$} & \multirow{2}*{$(E_{9},E_9)$}  & $(0;0,0,0,0,0,0,0,0,1,0,0,0,0,0,0,0,0;0)$ \\ 
&   &   &  &  &  & $(0;0,0,0,0,0,0,0,0,0,1,2,3,4,5,6,4,2;3)$ \\ 
\hline
\hline
\multirow{4}*{$\langle 2 \rangle$} & $A_{17}$  & $\mathrm{Bl}_{18}\bP^2$   & $\bP^2$             & $(-9,9)$ &  $(E_{18},0)$  & $(1;1,2,3,3,3,3,3,3,3,3,3,3,3,3,3,2,1;1)$ \\ 
\cline{2-7}
& $D_{16} A_1$   & $\mathrm{Bl}_{17}\bP^2$    & $\bP^1\times\bP^1$  & $(-8,8)$ & $(E_{17},A_1)$  & $(2;0,2,4,4,4,4,4,4,4,4,4,4,4,4,4,2,1;2)$ \\ 
\cline{2-7}
& $E_7D_{10}$    & $\mathrm{Bl}_{11}\bP^2$     & $\mathrm{Bl}_{7}\bP^2$            &  $(-2,2)$ & $(E_{11},E_7)$ & $(1;0,1,2,2,2,2,2,2,2,2,2,3,4,5,6,4,2;3)$ \\
\cline{2-7}
& $E_8E_8 A_1$  & $\mathrm{Bl}_{10}\bP^2$            & $\mathrm{Bl}_{8}\bP^2$ &  $(-1,1)$ & $(E_{10},E_8)$ & $(3;2,4,6,5,4,3,2,1,1,1,2,3,4,5,6,4,2;3)$ \\ 
\hline
\hline
\multirow{11}*{$\langle 4 \rangle$} & $A_{15}A_1A_1$ & $\mathrm{Bl}_{17}\bP^2$   & $\bP^1 \times \bP^1$             & $(-8,8)$ &  $(E_{17},A_1)$  & $(2;1,2,4,4,4,4,4,4,4,4,4,4,4,4,4,2,1;2)$  \\ 
\cline{2-7}
& $D_{17}$   & $\mathrm{Bl}_{18}\bP^2$    & $\bP^2$  & $(-9,9)$ & $(E_{18},0)$  & $(3;0,3,6,6,6,6,6,6,6,6,6,6,6,6,6,4,2;2)$ \\ 
\cline{2-7}
& $D_{16}\langle -4 \rangle$    & $\mathrm{Bl}_{17}\bP^2$    & $\bP^1\times\bP^1$  & $(-8,8)$ & $(E_{17},A_1)$  & $(3;0,3,6,6,6,6,6,6,6,6,6,6,6,6,6,3,1;3)$ \\ 
\cline{2-7}
& $D_{12}D_5$  & $\mathrm{Bl}_{13}\bP^2$            & $\mathrm{Bl}_{13}\bP^2$  &  $(-4,4)$ & $(E_{13},D_5)$ & $(2;0,2,4,4,4,4,4,4,4,4,4,4,4,5,6,4,2;3)$ \\ 
\cline{2-7}
 & $D_8D_8\langle -4 \rangle$  & $\mathrm{Bl}_{9}\bP^2$   &  $\mathrm{Bl}_{9}\bP^2$  & $(0,0)$ &  $(E_{9},E_9)$  & $(4;2,5,8,7,6,5,4,3,2,2,2,2,2,2,2,1,0;1)$ \\ 
\cline{2-7}
& \multirow{2}*{$E_8E_8\langle -4 \rangle$}   & $\mathrm{Bl}_{10}\bP^2$    & $\mathrm{Bl}_{8}\bP^2$  & $(-1,1)$ & $(E_{10},E_8)$  & $(9;6,12,18,15,12,9,6,3,1,0,0,0,0,0,0,0,0;0)$ \\ 
\cline{3-7}
&    & $\mathrm{Bl}_{9}\bP^2$     & $\mathrm{Bl}_{9}\bP^2$             &  $(0,0)$ & $(E_{9},E_9)$ & $(6;4,8,12,10,8,6,4,2,1,1,2,3,4,5,6,4,2;3)$ \\
\cline{2-7}
& $E_8D_9$  & $\mathrm{Bl}_{10}\bP^2$            & $\mathrm{Bl}_{8}\bP^2$ &  $(-1,1)$ & $(E_{10},E_8)$ & $(7;4,9,14,12,10,8,6,4,2,0,0,0,0,0,0,0,0;0)$ \\ 
\cline{2-7}
 & $E_7E_7A_3$  & $\mathrm{Bl}_{11}\bP^2$   & $\mathrm{Bl}_{7}\bP^2$             & $(-2,2)$ &  $(E_{11},E_7)$  & $(6;4,8,12,10,8,6,4,3,2,1,0,0,0,0,0,0,0;0)$ \\ 
\cline{2-7}
& \multirow{2}*{$E_6A_{11}$}   & $\mathrm{Bl}_{18}\bP^2$    & $\bP^2$  & $(-9,9)$ & $(E_{18},0)$  & $(3;2,4,6,5,4,3,3,3,3,3,3,3,3,3,3,2,1;1)$ \\ 
\cline{3-7}
&& $\mathrm{Bl}_{12}\bP^2$     & $\mathrm{Bl}_{6}\bP^2$            &  $(-3,3)$ & $(E_{12},E_6)$ & $(4;3,6,9,8,7,6,5,4,3,2,1,0,0,0,0,0,0;0)$ \\
\hline
\end{tabular}
\caption{Tyurin degenerations of K3 surfaces.}
\label{tab:degenerations}
\end{table}
\endgroup
\end{landscape}

\bibliography{publications,preprints}
\bibliographystyle{amsalpha}
\end{document}